# О нетривиальности быстрых (ускоренных) рандомизированных методов

А.В. Гасников, П.Е. Двуреченский, Усманова И.Н.

(ИППИ РАН, ПреМоЛаб ФУПМ МФТИ)

**Аннотация**

В данной работе предлагаются способы получения ускоренных и неускоренных вариантов рандомизированных покомпонентных методов и неускоренных вариантов методов рандомизации суммы, исходя из оптимальных методов для общих задач (стохастической) выпуклой оптимизации. В работе подчеркивается нетривиальность оценок, полученных для соответствующих ускоренных вариантов этих методов, которые выводятся в статье с помощью недавно предложенной техники каплинга. В отличие от многих других ситуаций, в данном случае не удается "вытащить", не погружаясь в детали доказательства (должным образом корректируя его), оптимальные методы (оценки) для рандомизированных покомпонентных методов и методов с рандомизацией суммы исходя из оптимальных методов (оценок), применимых к общим задачам стохастической оптимизации.

**Ключевые слова:** рандомизированные покомпонентные методы, быстрый градиентный метод, рандомизация суммы.

## 1. Введение

Для быстрого первоначального погружения в описываемую далее проблематику можно рекомендовать пп. 6.2 – 6.5 обзора [1].

В данной работе мы хотим подчеркнуть нетривиальность таких методов как, например, ускоренный (быстрый) покомпонентный метод Ю.Е. Нестерова [2], APPROX или ALPHA [3], [4], которые, в частности, являются покомпонентными вариантами быстрого градиентного метода (БГМ) [1] – эту ссылку можно также рекомендовать с точки зрения интересной подборки ссылок на работы, в которых объясняется, что такое БГМ. Нетривиальность в том, что этим методы являются рандомизированными и при этом ускоренными. Число необходимых итераций (как функция от желаемой точности) для таких ускоренных покомпонентных методов увеличивается в число раз $n$, равное размерности пространства, по сравнению с классическим БГМ, что и не удивительно, поскольку вместо всех $n \gg 1$ компонент градиента на каждой итерации используется только одна. Также нетривиальность в том, что если полный расчет градиента, скажем, требовал полного умножения разреженной матрицы на вектор – $sn$ операций, то пересчет (важно, что именно пересчет, а не расчет) компоненты градиента в определенных ситуациях требует всего $s$ операций (см. стр. 16–17 [3] и пп. 4, 5). Таким образом, увеличение числа итераций в $n$ раз компенсируется уменьшением стоимости одной итерации в $n$ раз (в не разреженном случае оговорка об "определенных ситуациях" существенно ослабляется, см. пп .4, 5). Но выгода от использования покомпонентных методов, как правило, есть из-за того, что в таких методах вместо константы Липшица градиента по худшему направлению (максимального собственного значения матрицы Гессе функционала) в оценки числа итераций входит "средняя" константа Липшица (оценивающаяся сверху средним арифметическим суммы диагональных элементов (следа) матрицы Гессе, т.е. средним арифметическим всех собственных чисел матрицы Гессе). Разница в этих константах для матриц Гессе, состоящих из элементов одного порядка, может равняться по порядку $n$ (см. пример 2 п. 5). На данный момент известно довольного много примеров применения покомпонентных методов для решения задач огромных размеров, в частности, приложений для задач моделирования сетей больших размеров и анализе данных [5].

В пп. 2, 3 мы демонстрируем те сложности, которые возникают при попытках получить ускоренные покомпонентные методы из оптимальных методов для задач стохастической оптимизации без погружения в то, как устроены эти оптимальные методы. Мы не ставили себе в этих пунктах цель: получить и подробно исследовать какие-то новые эффективные методы, поэтому изложение в этих пунктах ведется на "физическом" уровне строгости. В п. 4 мы приводим новое доказательство оценки скорости сходимости ускоренного покомпонентного метода, базирующееся на конструкции линейного каплинга [6]: БГМ = "выпуклая комбинация" прямого градиентного метода (ПГМ) и метода зеркального спуска (МЗС). Основная идея получения ускоренного покомпонентного метода: заменить в таком представлении БГМ в методах ПГМ и МЗС градиенты на соответствующие несмещенные оценки градиентов, полученные на основе покомпонентной рандомизации. Несмотря на то, что основной результат п. 4 (теорема 2 о сходимости предложенного метода и замечания к ней) – не есть полностью новый результат, подобные оценки (в различных частных случаях) ранее уже встречались в литературе, тем не менее, способ их получения (и его универсальность) представляется новым и весьма перспективным с точки зрения возможных последующих обобщений и приложений (некоторые примеры таких приложений и обобщений приведены в пп. 4, 5). Описанный способ также позволяет устанавливать различные новые факты об ускоренных покомпонентных методах. Наброски приведены в цикле замечаний к теореме 2 в п. 4 и в примерах п. 5. Однако этому планируется посвятить также и отдельную работу(-ы). В п. 6 кратко резюмируются результаты работы, приводятся заключительные замечания.

## 2. Нетривиальность ускоренных покомпонентных методов

Рассматривается задача гладкой выпуклой оптимизации

$$f(x) \to \min_{x \in Q}.$$

Мы постараемся сначала пояснить, как получить в неускоренном случае для этой задачи оценки для покомпонентных спусков из оценок рандомизированных методов решения этой задачи. Оказывается это можно довольно изящно сделать. К сожалению, при этом даже из оптимальных рандомизированных методов не удается "вытащить" оценки для ускоренных покомпонентных методов. В этом-то и заключается нетривиальность ускоренных покомпонентных методов.

Для простоты считаем, что везде в дальнейшем (в пп. 2, 3) мы говорим о 2-норме и евклидовой прокс-структуре (интересно было бы понять, как все, что далее будет написано, распространяется на более общие нормы/прокс-структуры [9]). Считаем, что функция $f(x)$ имеет Липшицев градиент с константой $L$, является $\mu$-сильно выпуклой, а множество $Q$ имеет диаметр $R$, при этом в точке минимума $\nabla f(x_*) = 0$. Последнее предположение – обременительное для задач условной оптимизации. К сожалению, мы пока не знаем как от него отказаться.

Будем считать, что на каждой итерации оракул выдает нам несмещенную оценку градиента с дисперсией $D$. Определим зависимость $N(\varepsilon)$ для изучаемого итерационного процесса: $N(\varepsilon)$ – наименьшее $N$ такое, что ($f_* = f(x_*)$, где $x_*$ – решение задачи)

$$E\left[f(x_N)\right] - f_* \le \varepsilon.$$

**Теорема 1.** *Существуют такие неускоренные методы (см., например, п. 6.2 [1]), которые работают по оценкам ($\Delta f^0 = f(x_0) - f_*$)*

$$\varepsilon = \min\left\{\mathrm{O}\left(\frac{LR^2}{N} + \sqrt{\frac{DR^2}{N}}\right), \mathrm{O}\left(\Delta f^0 \exp\left(-N\frac{\mu}{L}\right) + \frac{D}{\mu N}\right)\right\}.$$

*Существуют такие ускоренные методы (например, линейки SIGMA, см., например, [7]–[9]), которые работают по следующим не улучшаемым оценкам*

$$\varepsilon = \min\left\{\mathrm{O}\left(\frac{LR^2}{N^2} + \sqrt{\frac{DR^2}{N}}\right), \mathrm{O}\left(\Delta f^0 \exp\left(-N\sqrt{\frac{\mu}{L}}\right) + \frac{D}{\mu N}\right)\right\}.$$

Будем говорить о сильно выпуклом случае, если минимум в этих формулах достигается на втором аргументе.

Далее заметим, что если мы вместо обычного градиента используем его аппроксимации, возникающие в безградиентных и покомпонентных подходах [9] (когда оракул может на каждой итерации выдавать только значение функции в двух точках или производную по указанному нами направлению)

$$g_\tau(x,s) = \frac{n}{\tau}\left(f(x+\tau s) - f(x)\right)s \text{ или } g(x,s) = n\langle\nabla f(x), s\rangle s,$$

где $s$ – случайный вектор, равномерно распределенный на $S_2^n(1)$ – единичной сфере в 2-норме в пространстве $\mathbb{R}^n$, то имеет место следующий простой факт, являющийся следст-вием явления концентрации равномерной меры на сфере вокруг экватора [10]–[12] (северный полюс задается градиентом).

**Утверждение 1.** *Имеют место следующие формулы (см. [9]–[11])*

$$E_s\left[\left\|g_\tau(x,s)\right\|_2^2\right] \le 4n\left\|\nabla f(x)\right\|_2^2 + L^2\tau^2 n^2,$$

$$E_s\left[\left\|g(x,s)\right\|_2^2\right] = n\left\|\nabla f(x)\right\|_2^2.$$

Далее мы ограничимся оценкой для покомпонентного метода. "Почувствовать" эту оценку можно на еще более простом примере, когда несмещенная оценка градиента

$$g(x,s) = n\langle\nabla f(x), s\rangle s,$$

получается за счет другого выбора случайного вектора $s$. Мы считаем, что $s$ принимает равновероятно одно из $n$ направлений, соответствующих единичных ортов. В таком случае также имеет место соотношение

$$E_s\left[\left\|g(x,s)\right\|_2^2\right] = n\left\|\nabla f(x)\right\|_2^2.$$

Попробуем теперь исходя из полученных оценок и теоремы 1 в сильно выпуклом случае с $\nabla f(x_*) = 0$ получить оценку

$$N(\varepsilon) = \mathrm{O}\left(n\frac{L}{\mu}\ln\left(\frac{\Delta f^0}{\varepsilon}\right)\right)$$

для спусков по направлению. Итак, мы считаем, что вместо градиента оракул на каждой итерации (в точке $x_k$) может нам выдавать только производную по указанному нами направлению. Тогда мы имеем несмещенную оценку градиента с дисперсией (приводимая ниже оценка является не улучшаемой с точностью до мультипликативной константы)

$$D = \mathrm{O}\left(n\left\|\nabla f(x_k)\right\|_2^2\right).$$

Используя тот факт, что для любой гладкой выпуклой функции (в предположении, что $\nabla f(x_*) = 0$; для последнего неравенства еще нужно потребовать, чтобы $k$ не был слишком маленьким)

$$\left\|\nabla f(x_k)\right\|_2^2 \le 2L\cdot(f(x_k) - f_*) \le 2L\Delta f^0,$$

получаем из теоремы 1, что после

$$N = \mathrm{O}(5nL/\mu)$$

итераций

$$f(x_N) - f_* \leq \mathrm{O}\left(\Delta f^0 \exp\left(-N\frac{\mu}{L}\right) + \frac{2nL\Delta f^0}{\mu N}\right) \leq \mathrm{O}\left(\frac{1}{2}\Delta f^0\right)$$

и, тем более,

$$f(x_N) - f_* \leq \mathrm{O}\left(\Delta f^0 \exp\left(-N\sqrt{\frac{\mu}{L}}\right) + \frac{2nL\Delta f^0}{\mu N}\right) = \mathrm{O}\left(\frac{1}{2}\Delta f^0\right).$$

Тут также можно пользоваться методами, которые работают по оценкам [13]–[15]

$$\varepsilon = \mathrm{O}\left(M^2/(\mu N)\right),$$

где

$$M^2 = E_s\left[\left\|g_{\tau,\delta}(x_k, s)\right\|_2^2\right] = \mathrm{O}\left(n\left\|\nabla f(x_k)\right\|_2^2\right),$$

что не удивительно, поскольку мы, фактически, при данном подходе и работаем с $M^2$, а не с $D \leq M^2$. В любом случае при таком способе рассуждений возникает неаккуратность, связанная с тем, что мы лишь обеспечили

$$f(x_N) - f_* \leq \mathrm{O}\left(\frac{1}{2}\Delta f^0\right).$$

В действительности, тут нужно аккуратно выписывать константы, которые в итоге увеличат константу "5" в ожидаемой сейчас формуле $N = 5nL/\mu$ в несколько раз. Наконец, необходимо проводить рассуждения с оценками вероятностей больших уклонений (здесь помогают самые грубые неравенства типа Буля, поскольку, имеются субгауссовские хвосты, а точнее вообще финитный носитель у стохастического градиента). Далее мы уже не будем делать такие оговорки, поскольку в этом и следующем пункте мы преследуем цель – продемонстрировать нетривиальность ускоренных рандомизированных методов, а не точного выписывания методов, которые получаются по ходу рассуждений. Эти методы не очень интересны, поскольку заведомо не являются оптимальными.

Делая $\log_2(\Delta f^0/\varepsilon)$ таких перезапусков (стартуем в новом цикле с той точки, на которой остановились на прошлом цикле) с $N = \mathrm{O}(5nL/\mu)$ итерациями на каждом перезапуске (цикле), в итоге получим оценку общего числа итераций

$$N(\varepsilon) = \mathrm{O}\left(n\frac{L}{\mu}\ln\left(\frac{\Delta f^0}{\varepsilon}\right)\right).$$

Здесь можно "поиграться" на, так называемом, mini-batch'инге (см., например, п. 6.2 [1]), для этого нужна уже формула с дисперсией, то есть $\varepsilon = \mathrm{O}\left(M^2/(\mu N)\right)$ не подходит.

Тем не менее, даже при использовании не улучшаемых (с точностью до мультипликативной константы) рандомизированных методов из теоремы 1 мы не смогли получить оценки работ [2]–[4], для ускоренных покомпонентных методов

$$N(\varepsilon) = \mathrm{O}\left(n\sqrt{\frac{L}{\mu}}\ln\left(\frac{\Delta f^0}{\varepsilon}\right)\right).$$

Строго говоря, в работах [3], [4] таких оценок для покомпонентных методов (в сильно выпуклом случае) мы и не видели (отметим при этом, что такая оценка есть для метода из работы [2] и для безградиентного метода из работы [16]), однако в [3], [4] есть аналогичные "ускоренные" оценки в несильно выпуклом случае. Кроме того, Питер Рихтарик [17] сообщил нам, что он умеет устанавливать эти оценки (для ускоренных покомпонентных методов и в сильно выпуклом случае), и сейчас готовит статью на эту тему (свой способ получения таких оценок мы изложим в п. 4). Проблема тут в том, что мы пользовались правым неравенством (считаем $\nabla f(x_*) = 0$)

$$2\mu \cdot \left(f(x_k) - f_*\right) \leq \left\|\nabla f(x_k)\right\|_2^2 \leq 2L \cdot \left(f(x_k) - f_*\right).$$

Мы специально здесь написали и левое неравенство. Отсюда видно, что в принципе при использовании правого неравенства мы можем потерять $L/\mu$. Не удивительно, что в итоге мы, действительно, теряем $\sqrt{L/\mu} \ll L/\mu$. К сожалению, такого рода рассуждения не позволяют никак "вытащить" оценки оптимального (ускоренного) покомпонентного метода из соответствующих оптимальных полных градиентных методов, не погружаясь в детальный анализ доказательства их сходимости. Нетривиально то, что это один из тех редких примеров (другой см. ниже в п. 3), когда такая философия переноса не сработала. Обычно все удается перенести без особых погружений в детали доказательства. То есть работает принцип: *оптимальный метод порождает оптимальный*.[1]

### 3. Нетривиальность ускоренных методов рандомизации суммы

Рассмотрим теперь в задаче п. 2 случай когда (этот случай разбирается, например, в п. 6.3 [1] и работах [18]–[21])

---

[1] Хорошо известные примеры тут: 1) регуляризация [2], позволяющая переносить оптимальные методы, работающие в сильно выпуклом случае, на просто выпуклый случай; 2) техника рестартов [6], [14] (см. также п. 4 далее), позволяющая из оптимального метода для выпуклой задачи, получить оптимальный метод для сильно выпуклой задачи.

$$f(x) = \frac{1}{m}\sum_{k=1}^{m} f_k(x),$$

где все функции гладкие с константой Липшица градиента $L$. Также как и раньше считаем $f(x)$ $\mu$-сильно выпуклой. В качестве несмещенной оценки градиента будем брать вектор (по поводу определения SIGMA см. [8])

$$\nabla f(x_t^s, \xi) = \nabla f_\xi(x_t^s) - \nabla f_\xi(y^s) + \nabla f(y^s), \ y^s = x_N^{s-1},$$

$$x_{t+1}^s = \text{SIGMA}(x_t^s, \nabla f(x_t^s, \xi)), \ t = 0,\ldots,N-1,$$

где случайная величина $\xi$ принимает равновероятно одно из значений $1,\ldots,m$; параметр $N$ будет выбран позже как $N = \mathrm{O}(4L/\mu)$. Здесь по $t$ идет внутренний цикл, а по $s$ внешний.

Приведенный метод (как и метод из п. 6.3 [1]) можно обобщить (с сохранением всех последующих оценок и способов их получения) на стохастический случай, когда $f_k(x) \coloneqq E_\eta[f_k(x;\eta_k)]$, где $f_k(x;\eta_k)$ – выпуклые по $x$ функции с равномерно (по $\eta_k$) ограниченными (числом $L$) константами Липшица градиентов. При этом (случайная величина $\eta_\xi^{s,t}$ имеет такое же распределение, как и $\eta_\xi$; также считаем, что $\eta_\xi^{s,t}$ ни от чего не зависит, в частности, от других $\{\eta_\xi^{s,t}\}$ и от $\xi$)

$$\nabla f(x_t^s, \xi; \eta_\xi^{s,t}) = \nabla f_\xi(x_t^s; \eta_\xi^{s,t}) - \nabla f_\xi(y^s; \eta_\xi^{s,t}) + \nabla f(y^s).$$

Это очевидно, для случая, когда

$$f_k(x) \coloneqq E_\eta[f_k(x;\eta_k)] = \frac{1}{l}\sum_{i=1}^{l} f_k(x;i),$$

поскольку все сводится к исходной постановке с $m \coloneqq ml$. Далее мы ограничимся рассмотрением только детерминированного случая.

**Утверждение 2.** *Имеет место следующая оценка ($\Delta f^s = f(y^s) - f_*$)*

$$D^s = E_\xi\left[\left\|\nabla f(x_t^s, \xi) - E_\xi[\nabla f(x_t^s, \xi)]\right\|_*^2\right] = \mathrm{O}\big(L\cdot(f(y^s) - f_*) + L\cdot(f(x_t^s) - f_*)\big) = \mathrm{O}(L\Delta f^s).$$

Для доказательства утверждения 2 в случае $\nabla f(x_*) = 0$ см., например, формулу (6.3) и лемму 6.1 [1] и цитированную в п. 6.3 [1] литературу (в частности, [18]–[20]). Причем это утверждение можно формулировать с точными константами вместо O( ), чтобы ей можно было далее практически воспользоваться (однако, мы здесь не будем этого де-лать). В случае $\nabla f(x_*) \neq 0$ доказательство утверждения 2 нам не известно (не известно даже останется ли оно верным).

Возьмем теперь в теореме 1 $N = \mathrm{O}(4L/\mu)$ и воспользуемся утверждением 2

$$\Delta f^{s+1} = \mathrm{O}\left(\Delta f^s \exp\left(-N\sqrt{\frac{\mu}{L}}\right) + \frac{D^s}{\mu N}\right),$$

$$\frac{D^s}{\mu N} = \mathrm{O}\left(\frac{L \Delta f^s}{\mu N}\right) = \mathrm{O}\left(\frac{1}{4}\Delta f^s\right), \quad \Delta f^s \exp\left(-N\sqrt{\frac{\mu}{L}}\right) = \mathrm{O}\left(\frac{1}{4}\Delta f^s\right).$$

Получим

$$\Delta f^{s+1} \leq \mathrm{O}\left(\frac{1}{2}\Delta f^s\right).$$

Здесь можно "поиграться" на mini-batch'инге [21], вычисляя вместо $\nabla f(x_t^s, \xi)$ агрегат

$$\frac{1}{r}\sum_{i=1}^{r} \nabla f(x_t^s, \xi_i).$$

Таким образом, если у нас есть возможность параллельно на одной итерации вычислять градиенты $\nabla f_{\xi_i}(x_t^s)$, то (поскольку дисперсия этого агрегата будет $\mathrm{O}(L\Delta f^s/r)$) можно выбирать $r = \mathrm{O}(2\sqrt{L/\mu})$ (выбрать $r$ большим нельзя) и, соответственно, сократить число итераций на цикле до $N = \mathrm{O}(2\sqrt{L/\mu})$. Далее, также как и раньше, делаем $\log_2(\Delta f^0/\varepsilon)$ перезапусков (циклов), на каждом из которых вначале надо посчитать полный градиент (это стоит $m$ вычислений градиентов слагаемых), а потом еще сделать $N$ итераций, на каждой из которых дополнительно требуется рассчитывать в новой точке градиент одного (или нескольких, если используется mini-batch'инг) слагаемого. Таким образом, общая сложность (измеряемая на этот раз не в итерациях, а числе вычислений градиентов слагаемых в сумме в представлении $f(x)$, при этом мы считаем, что сложность вычисления разных слагаемых одинакова в категориях $\mathrm{O}(\ )$) будет

$$\mathrm{O}\left(\left(m + \frac{L}{\mu}\right) \cdot \ln\left(\frac{\Delta f^0}{\varepsilon}\right)\right),$$

что соответствует части нижней оценки в классе детерминированных методов (см. [22])

$$\mathrm{O}\left(\left(m + \min\left\{\sqrt{m\frac{L}{\mu}}, \frac{L}{\mu}\right\}\right) \cdot \ln\left(\frac{\Delta f^0}{\varepsilon}\right)\right).$$

Шалевым-Шварцем и Зангом был поставлен вопрос (см., например, п. 6.3 [1]): возможно ли достичь такой (нижней) оценки в целом каким-нибудь методом? То есть речь опять (как и раньше) на самом деле о том, можно ли сохранить ускоренность метода? Только пока, наверное, не очень понятно, причем здесь (ускоренные) рандомизированные покомпонентные методы. Ситуация проясняется, если мы выделим строго выпуклое слагаемое из $f(x)$ в отдельный композит $g(x)$ (см., например, п. 5.1 [1]):

$$f(x) = \frac{1}{m}\sum_{k=1}^{m} f_k(x) + g(x),$$

и построим специальную двойственную задачу (см. замечание 5 [5], а также пример 4 п. 5 ниже). Тогда (при некоторых дополнительных предположениях [5], [23], [24]) удается показать, что в некотором смысле выписанная нижняя оценка, действительно, достигается на ускоренном покомпонентном методе для двойственной задачи. Только для сопоставления потребуется перейти от анализа числа вычислений слагаемых к общему числу арифметических операций (см. пример 4 п. 5 ниже). Собственно, Шалев-Шварц и Занг в работе [23] сами таким образом (правда, не подчеркивая, что, по сути, используют для двойственной задачи ускоренный покомпонентный спуск) и привели пример достижимости (опять оговоримся, что в некотором смысле) нижней оценки.

В связи со сказанным выше отметим, что приведенные в этом разделе оценки (подобно оценкам для покомпонентных методов, см., например, п. 4) можно рассматривать в случае не равномерной рандомизации (выбора слагаемых), а также в случае разных свойств гладкости у разных функций. Сейчас в оценки методов, описанных в этом пункте, входит худшая (по всем слагаемым) константа Липшица градиентов (мы выбрали все константы Липшица градиентов слагаемых для наглядности одинаковыми). В действительности, можно перейти к некоторым их средним вариантам. Однако мы не будем здесь этого делать, поскольку, как уже отмечалось, природа описанных рандомизированных методов вскрывается применением покомпонентных методов к двойственной задачи, для которых все эти нюансы хорошо проработаны.

В этом примере мы опять видим, что попытка из оптимального метода для задач стохастической оптимизации SIGMA (см. теорему 1) вытащить оптимальные оценки для задачи со специальной (скрытой, через двойственную задачу, покомпонентной) структурой не привели к успеху (на mini-batch'инг не стоит обращать внимание, он просто позволяет параллелить вычисления, не более того). Таким образом, это лишний раз подчеркивает некоторую самостоятельность и важность отдельного изучения ускоренных покомпонентных методов. Именно такого типа методы (например, APPROX, ALPHA [3], [4]) позволяют получать наилучшие оценки. И оценки скорости сходимости этих методов представляют отдельный интерес (в смысле их получения). Насколько нам известно (лучше всего следить за этой областью по работам П. Рихтарика [17] и Т. Занга [25], и цитированной ими литературы), сейчас для таких методов используется только евклидова прокс-структура, используются только простые ограничения (сепарабельные), которые обычно зашивают в композитный член [3]. Интересно также было бы понять (охарактеризовать) класс задач, в которых возможно эффективно организовать пересчет компоненты градиента для ускоренных покомпонентных методов. Кое-что на эту тему имеется вот здесь стр. 16–17 [3], [26], слайды 6–9 [27]. Подробнее все это будет рассмотрено далее.

## 4. Получение ускоренных покомпонентных методов с помощью каплинга неускоренных прямых покомпонентных методов и покомпонентного метода зеркального спуска

Как уже отмечалось на текущий момент не до конца ясно, насколько все, что сейчас известно для методов, в которых доступен полный градиент, имеет свои аналоги и в (блочно-)покомпонентных методах. Скажем, не все понятно с тем, как можно играть на выборе прокс-функции в ускоренных (блочно-)покомпонентных методах, не до конца ясно: можно ли (если можно, то каким образом) использовать ускоренные покомпонентные методы, если рассматривается задача условной минимизации с множеством специальной простой структуры (в смысле прокс-проектирования, а точнее (блочно-)покомпонентной версии этой операции), но, вообще говоря, не сепарабельной структуры – в частности, в таких задачах $\nabla f(x_*) \neq 0$ [2], [3]; есть ли аналог универсального метода Ю.Е. Нестерова [28] в покомпонентном варианте; имеют ли покомпонентные методы прямо-двойственную структуру [5]; как перенести на покомпонентные методы концепцию неточного оракула Деволдера–Глинера–Нестерова (см., например, [7]–[9]); верно ли, что для ускоренного покомпонентного метода расстояние от любо точки итерационного процесса до решения всегда ограничено некоторой универсальной небольшой константой (меньшей 10), умноженной на расстояние от точки старта до решения, как это имеет место для БГМ (замечание 4 [9])? Список можно продолжить, однако мы здесь остановимся, и сформулируем общий тезис: *все, что сейчас известно для методов первого порядка, в которых доступен полный градиент, имеет (с оговорками, о возможности перенесения результатов на не сепарабельные множества) свои аналоги и в (блочно-)покомпонентных методах; более того, константы (Липшица градиента), фигурирующие в обычных градиентных методах, рассчитанные на худший случай (худшее направление), в покомпонентных методах заменяются "средними" значениями, что в определенных ситуациях может давать ускорение в корень из размерности пространства раз (не говоря о том, что покомпонентные методы при этом могут еще и хорошо параллелиться)* [3], [5], [27], [29].

К сожалению, исходя из всех известных нам опубликованных на данный момент способов вывода (доказательства сходимости) покомпонентных методов (наиболее, конечно, интересны тут ускоренные варианты) сформулированный выше тезис (гипотеза) совсем не кажется очевидным. Однако совсем недавно, в работе [6] был предложен изящный и перспективный во многих отношениях[2] способ получения БГМ с помощью выпук-лого каплинга (комбинации) обычного (неускоренного) прямого градиентного метода[3] (ПГМ) и метода зеркального спуска (МЗС). Ряд "хороших" свойств (например, прямо-двойственность) "наследуются" при таком представлении от зеркального спуска. Естественно, возникает идея попробовать использовать соответствующие легко исследуемые в отдельности покомпонентные аналоги этих двух структурных блоков, чтобы получить ускоренный покомпонентный метод. Оказывается, что это, действительно, можно сделать. Далее, основываясь на результатах работы [6], мы приведем соответствующие выкладки.

---

[2] В том числе в отношении более простого обоснования отмеченной выше возможности перенесения свойств с полноградиентных методов на покомпонентные.

Исходя из написанного в предыдущих пунктах, можно сказать, что для получения ускоренных покомпонентных методов требуется более тонкая игра (на каждой итерации) на правильном сочетании базовых методов со специальным выбором параметров. Оптимальный метод порождается выпуклой комбинацией неоптимальных методов для класса гладких задач, и именно из этого стоит исходить (распространяя конструкцию на покомпонентные методы), чтобы получить ускоренный покомпонентный метод. Этот тезис нам также представляется полезным, поскольку он подтверждает, что оптимальные методы порождают оптимальные, просто в ряде случаев требуется заглядывание в структуру (базис) метода, чтобы иметь возможность из него породить что-то новое оптимальное.

Сначала мы постараемся в максимально упрощенной ситуации пояснить, как можно получить ускоренный покомпонентный метод, исходя из конструкции п. 3 работы [6]. Все, что далее будет написано, допускает серьезные обобщения, о которых мы упомянем ближе к концу этого пункта.

Итак, рассмотрим задачу

$$f(x) \to \min_{x \in Q}.$$

Введем необходимые в дальнейшем обозначения/определения:

$$e_i = \overbrace{(0,...,0,\underbrace{1}_{i},0,...,0)}^{n};$$

$\left| \partial f(x + he_i)/\partial x_i - \partial f(x)/\partial x_i \right| \le L_i h$ для всех $x \in \mathbb{R}^n$ и $h \in \mathbb{R}$;

$$\|x\|^2 = \sum_{i=1}^{n} L_i x_i^2, \quad \|\nabla f(x)\|_*^2 = \sum_{i=1}^{n} L_i^{-1} \left(\partial f(x)/\partial x_i\right)^2;$$

---

[3] Использование в БГМ в качестве одного из структурных блоков именно ПГМ (это заметно упрощает рассуждения в случае $Q = \mathbb{R}^n$ по сравнению с другими возможными вариантами) не является обязательным атрибутом. По-видимому, можно построить (в схожем ключе) аналог БГМ (с аналогичными оценками скорости сходимости) на базе МЗС (или его "сходящегося" варианта [30]) и прямого проксимального градиентного метода (ППГМ) и(или) двойственного градиентного метода [7], который будет лишен отмеченных недостатков.

$$d(x) = \frac{1}{2}\|x\|^2, \ V_x(y) = d(y) - \langle \nabla d(x), y-x \rangle - d(x) = \frac{1}{2}\|y-x\|^2;$$

$$\nabla_i f(x) = \underbrace{(0,...,0, \partial f(x)/\partial x_i, 0,...,0)}_{i};$$

$i \in [1,...,n]$ – означает, что $P(i = j) = n^{-1}, \ j = 1,...,n$;

$E_i[G(i)]$ – математическое ожидание по $i \in [1,...,n]$;

$E_{i_{k+1}}[G(i_1,...,i_{k+1})|i_1,...,i_k] = g(i_1,...,i_k)$ – условное математическое ожидание по $i_{k+1} \in [1,...,n]$;

$E_{i_1,...,i_k}\left[E_{i_{k+1}}[G(i_1,...,i_{k+1})|i_1,...,i_k]\right] = E[G(i_1,...,i_{k+1})]$ – полное математическое ожидание по всему набору $i_1,...,i_{k+1} \in [1,...,n]$;

$$\text{Grad}_i(x) = \arg\min_{\tilde{x} \in Q}\left\{\langle \nabla_i f(x), \tilde{x} - x \rangle + \frac{1}{2}\|\tilde{x} - x\|^2\right\} = x - \frac{1}{L_i}\nabla_i f(x);$$

$$\text{Mirr}_z(\xi) = \arg\min_{y \in Q}\{\langle \xi, y - z \rangle + V_z(y)\} = \left(\left\{z_i - \frac{1}{L_i}\xi_i\right\}_{i=1}^n\right).$$

Приведенные формулы специально были записаны таким образом, чтобы их легко было перенести на случай, когда выбирается не одна компонента $i$, а целый блок компонент и $Q \neq \mathbb{R}^n$.

Опишем костяк покомпонетного ускоренного метода (Accelerated by Coupling Randomized Coordinate Descent – ACRCD) на базе специального каплинга покомпонентных вариантов ПГМ (Grad) и МЗС (Mirr) ($x_0 = y_0 = z_0$)

$$\text{ACRCD}(\alpha, \tau; \Theta, x_0, f(x_0) - f_* \leq d)$$

1. $x_{k+1} = \tau z_k + (1-\tau) y_k$, $\tau \in [0,1]$ – будет выбрано позже;

2. $i_{k+1} \in [1,...,n]$ – независимо от предыдущих розыгрышей;

3. $y_{k+1} = \text{Grad}_{i_{k+1}}(x_{k+1})$;

4. $z_{k+1} = \text{Mirr}_{z_k}(\alpha n \nabla_{i_{k+1}} f(x_{k+1}))$, $\alpha > 0$ – будет выбрано позже.

Поскольку

$$E_i[n \nabla_i f(x)] = \nabla f(x),$$

то шаг 4 (согласно формуле (3.1) [6]; в евклидовом случае можно ограничиться более простыми рассуждениями – см., например, стр. 223 [31]; другой способ получить оценки для МЗС, показывающий дополнительную связь МЗС и проксимального ПГМ (ППГМ), – воспользоваться оценками для метода ППГМ с неточным оракулом из [7], [28]) влечет[4]

$$\alpha n \left\langle \nabla_{i_{k+1}} f(x_{k+1}), z_k - u \right\rangle \leq \frac{\alpha^2 n^2}{2} \left\| \nabla_{i_{k+1}} f(x_{k+1}) \right\|_*^2 + V_{z_k}(u) - V_{z_{k+1}}(u) =$$

$$= \frac{\alpha^2 n^2}{2 L_{i_{k+1}}} \left| \partial f(x_{k+1}) / \partial x_{i_{k+1}} \right|^2 + V_{z_k}(u) - V_{z_{k+1}}(u) \overset{\text{шаг 3}}{\leq} \alpha^2 n^2 \left( f(x_{k+1}) - f(y_{k+1}) \right) + V_{z_k}(u) - V_{z_{k+1}}(u).$$

Возьмём от этого неравенства условное математическое ожидание $E_{i_{k+1}}\left[ \cdot \,|\, i_1, ..., i_k \right]$:

$$\alpha \left\langle \nabla f(x_{k+1}), z_k - u \right\rangle \leq \alpha^2 n^2 \left( f(x_{k+1}) - E_{i_{k+1}}\left[ f(y_{k+1}) \,|\, i_1, ..., i_k \right] \right) + V_{z_k}(u) - E_{i_{k+1}}\left[ V_{z_{k+1}}(u) \,|\, i_1, ..., i_k \right].$$

Согласно формуле (3.2) [6], которая используется в совершенно таком же виде, как и в [6], из последнего неравенства при

$$\frac{1-\tau}{\tau} = \alpha n^2,$$

получаем

$$\alpha \left\langle \nabla f(x_{k+1}), x_{k+1} - u \right\rangle \leq \alpha^2 n^2 \left( f(y_k) - E_{i_{k+1}}\left[ f(y_{k+1}) \,|\, i_1, ..., i_k \right] \right) + V_{z_k}(u) - E_{i_{k+1}}\left[ V_{z_{k+1}}(u) \,|\, i_1, ..., i_k \right].$$

Положим $u = x_*$, и возьмём математическое ожидание $E_{i_1, ..., i_k}\left[ \cdot \right]$ (если $k \geq 1$) от каждого такого неравенства, и просуммируем то что получается по $k = 0, ..., K-1$[5]

$$\alpha K \left( E\left[ f(\bar{x}_K) \right] - f_* \right) \leq \alpha \sum_{k=1}^{K} E\left[ \left\langle \nabla f(x_k), x_k - x_* \right\rangle \right] \leq$$

$$\leq \alpha^2 n^2 \left( f(x_0) - E\left[ f(y_K) \right] \right) + V_{x_0}(x_*) - V_{z_K}(x_*) \leq \alpha^2 n^2 \left( f(x_0) - f_* \right) + V_{x_0}(x_*),$$

---

[4] Отметим, что первое неравенство специально записано таким образом (в достаточно общем виде) чтобы была видна возможность рассмотрение прокс-структур отличных от евклидовой. По-видимому, для покомпонентных методов в подавляющем большинстве приложений можно ограничиться рассмотрением только евклидовой прокс-структуры. Нюансы могут возникать, когда вместо одной компоненты разрешается сразу случайно выбирать целый блок компонент (необязательно постоянного размера) [2]–[4], что для евклидовой прокс-структуры можно также понимать как mini-batch'инг [1]. Иногда в приложениях это (использовать сразу блок случайно выбранных компонент, причем в понятие "случайно" тут можно много что вкладывать) бывает полезно [32].

[5] Из этого неравенства устанавливается прямо-двойственная природа метода ACRCD [33]. Прямо-двойственность ускоренных покомпонентных методов требуется в ряде приложений, см., например, [5]. Впрочем, в сильно выпуклом случае (к которому все можно сводить за дополнительную логарифмическую плату, см. далее) прямо-двойственность оказывается уже не нужна (см., например, главу 3 [7]).

$$\overline{x}_K = \overline{x}_K(x_0) = \frac{1}{K}\sum_{k=1}^{K} x_k.$$

Пусть

$$V_{x_0}(x_*) \le \Theta,$$

$$f(x_0) - f_* \le d.$$

Выбирая

$$\tau = \frac{1}{\alpha n^2 + 1}, \; \alpha = \frac{1}{n}\sqrt{\frac{\Theta}{d}}, \; K = K(d) = 8n\sqrt{\frac{\Theta}{d}}$$

получим

$$E\left[f(\overline{x}_K)\right] - f_* \le \frac{2n\sqrt{\Theta d}}{K} \le \frac{d}{4}.$$

Для получения сходимости по вероятности, воспользуемся следующим приемом [9], о котором мы узнали от А.С. Немировского. Из

$$E\left[f(\overline{x}_K)\right] - f_* \le d/4$$

по неравенству Маркова

$$X = f(\overline{x}_K) - f_* \ge 0, \; P(X \ge t) \le E[X]/t, \; t = d/2,$$

имеем

$$P\left(f(\overline{x}_K) - f_* \ge d/2\right) \le 1/2.$$

Отсюда следует, что если мы независимо (можно параллельно) запустим $\lceil \log_2(\sigma^{-1}) \rceil$ траекторий $\mathrm{ACRCD}(\alpha, \tau; \Theta, x_0, d)$ (определив тот $\overline{x}_K$, для которого значение $f(\overline{x}_K)$ будет наименьшим), то с вероятностью $\ge 1 - \sigma$ хотя бы на одной траектории будем иметь

$$f(\overline{x}_K) - f_* \le d/2.$$

К сожалению, это требует расчета в $\lceil \log_2(\sigma^{-1}) \rceil$ точках значения функции $f(x)$. Расчет функции в точке $f(x)$ может быть заметно дороже стоимости одной итерации метода ACRCD. Однако это все равно не изменит по порядку оценку общего числа арифметических операций.

Итак, пусть $\text{ACRCD}(\alpha,\tau;\Theta,x_0,d)$ выдал такое $\bar{x}_{K(d)}(x_0)$, что с вероятностью $\geq 1-\sigma$ имеет место неравенство

$$f\left(\bar{x}_K(x_0)\right) - f_* \leq d/2.$$

Важно заметить, что при этом с вероятностью $\geq 1-\sigma$

$$V_{\bar{x}_{K(d)}(x_0)}(x_*) \leq \max_{x \in Q}\left\{V_x(x_*):\ f(x)-f_* \leq d/2\right\} \leq \max_{x \in Q}\left\{V_x(x_*):\ f(x)-f_* \leq d\right\}.$$

При этом выписанная оценка не запрещает, например, что $V_{\bar{x}_{K(d)}(x_0)}(x_*) \gg \Theta$. В этой связи для правильной работы описываемой далее процедуры перезапусков, к сожалению, необходимо переопределить $\Theta$ следующим образом (все это может существенно ухудшить итоговую оценку,[6] однако впоследствии с помощью регуляризации исходной постановки задачи мы покажем, как можно практически полностью нивелировать эту проблему)

$$\max_{x \in Q}\left\{V_x(x_*):\ f(x)-f_* \leq d\right\} \leq \Theta.$$

Запустим далее $\text{ACRCD}(\alpha,\tau;\Theta,\bar{x}_K(x_0),d/2)$, получим такой $\bar{x}_{K(d/2)}(\bar{x}_K(x_0))$, что с вероятностью $\geq 1-2\sigma$ (воспользовались неравенством Буля) имеет место неравенство

$$f\left(\bar{x}_{K(d/2)}(\bar{x}_K(x_0))\right) - f_* \leq d/4.$$

Процесс можно продолжать ... Для достижения по функции точности $\varepsilon$ с вероятностью

$$\geq 1 - \lceil \log_2(d/\varepsilon) \rceil \sigma$$

будет достаточно $\lceil \log_2(d/\varepsilon) \rceil$ таких итераций (перезапусков). Требуемое при этом общее число итераций (число обращений к компонентам вектора градиента) оценивается сверху следующим образом:[7]

$$N \leq 8n\sqrt{\frac{\Theta}{\varepsilon}}\left(1 + 2^{-1/2} + 2^{-1} + 2^{-3/2} + \ldots\right)\log_2(\sigma^{-1}) < 30n\sqrt{\frac{\Theta}{\varepsilon}}\log_2(\sigma^{-1}).$$

---

[6] Например, для выпуклой квадратичной функции (безобидной, с первого взгляда, в виду равномерной ограниченности всех коэффициентов), о которой нам сообщил Ю.Е. Нестеров: $f(x) = x_1^2 + \sum_{k=1}^{n-1}(x_{k+1}-2x_k)^2$ множества Лебега оказываются сильно сплющенными (плохо обусловленными). Также эту функцию интересно прооптимизировать с помощью ПГМ методов с не евклидовой нормой [6]. Эти методы релаксационные (то есть значение функции монотонно убывает на итерациях), но при этом точки, генерируемые методами, по ходу итерационного процесса могут уходить намного дальше от решения, чем точка старта.

[7] Выбирая $K = 9n\sqrt{\Theta/d}$, можно уменьшить константу 30 до 27, последняя константа уже не улучшаемая при таком способе рассуждений.

**Теорема 2.** *После $\lceil \log_2(d/\varepsilon) \rceil$ описанных выше рестартов, метод ACRCD выдает такой $x^N$, что с вероятностью $\geq 1 - \sigma$ имеет место неравенство*

$$f(x^N) - f_* \leq \varepsilon.$$

*При этом методу требуется для этого сделать*

$$N = N(\varepsilon) = 27n\sqrt{\frac{\Theta}{\varepsilon}} \log_2\left(\frac{\log_2(d/\varepsilon)}{\sigma}\right)$$

*итераций.*

**Замечание 1. (достоинства и недостатки ACRCD)** Из описанной конструкции ACRCD, как уже отмечалось ранее, следует его прямо-двойственность [33]. Заметив, что ПГМ и МЗС (Grad и Mirr) легко могут быть обобщены на композитные постановки задач [34] (к этому случаю, в частности, можно свести и минимизацию на параллелепипеде), с помощью описанной выше конструкции можно получить соответствующий композитный вариант ACRCD для задач с сепарабельным композитом (см. также [3]). К сожалению, ряд других свойств ACRCD уже не так просто "вытащить" (и не всегда понятно даже, возможно ли это в принципе, и имеют ли вообще нужные свойства место здесь). В частности, например, не понятно, как можно адаптивно (по ходу итерационного процесса) подбирать константы Липшица по разным направлениям, подобно п. 6.1 работы [2]. Не понятно, как можно "бороться" с проблемой неизвестности одновременно двух параметров $\Theta$ и $d$, нужных методу для работы.[8] Описанный нами вариант метод ACRCD работает в предположении $Q = \mathbb{R}^n$, и не гарантирует свойство равномерной ограниченности (в вероятностях категориях [9], [32]) последовательности расстояний от решения до точек генерируемых методом, значением этого расстояния в начальный момент, умноженным не небольшую универсальную константу (в частности, не зависящую от свойств функционала задачи). Эти плохие свойства ACRCD[9] "унаследовал" от БГМ, описанного в п. 3 [6]. В п. 4 [6] описан вариант немного другой вариант БГМ, который лишен этих недостатков. Оказывается можно распространить и его на покомпонентный случай, что далее будет сделано (см. замечание 2).

Уже отмеченные свойства (прямо-двойственность, обобщение на композитные задачи) и все далее изложенные свойства (обобщения) ACRCD допускают всевозможные сочетания друг с другом. Детали мы вынуждены опустить (планируется посвятить этому отдельную работу), но, в большинстве случаев, все это является довольно простыми факта-ми (впрочем, как правило, требующие для аккуратного доказательства довольно громоздких, но вполне стандартных рассуждений). Можно сказать по-другому: далее приводится "базис" для всевозможных последующих обобщений.

---

[8] Стандартные приемы рестартартов по неизвестному параметру разработаны сейчас только для случая одного неизвестного параметра [9], формальная попытка перенесения на случай двух и более неизвестных параметров (без дополнительных предположений [32]) приводит к резкому увеличению сложности процедуры.

[9] Вместе с уже отмеченной ранее проблемой вхождения $\Theta$ в итоговую оценку (числа итераций $N(\varepsilon)$) вместо $V_{x_0}(x_*)$, как это можно было ожидать [3].

**Замечание 2. (ACRCD\*)** Используемая при построении ACRCD техника рестартов позволила довольно просто получить оценки вероятностей больших уклонений. Однако эта же техника создала ряд проблем (см. замечание 1), многие из которых, в первую очередь, связаны с некоторым запаздыванием в обновлении параметров $\tau$ и $\alpha$. Они обновляются только на новом рестарте. Основная идея (см. п. 4 [6]) – сделать эти параметры зависящими от шага. Тогда удастся избавиться от рестартов и приобрести ряд хороших свойств. Далее описывается соответствующая модификация метода ACRCD. Предварительно определим две числовые последовательности

$$\alpha_1 = \frac{1}{n^2}, \; \alpha_k^2 n^2 = \alpha_{k+1} n^2 - \alpha_{k+1}, \; \tau_k = \frac{1}{\alpha_{k+1} n^2}.$$

Можно написать явные формулы. Также можно, следуя [6], брать близкие последовательности (теоретически немного проще исследовать первый вариант, но второй вариант более нагляден, а итоговые оценки скорости сходимости практически идентичны)

$$\boxed{\alpha_{k+1} = \frac{k+2}{2n^2}, \; \tau_k = \frac{1}{\alpha_{k+1} n^2} = \frac{2}{k+2}.}$$

ACRCD\*$(x_0 = y_0 = z_0)$

$$\boxed{\begin{array}{l} 1. \; x_{k+1} = \tau_k z_k + (1-\tau_k) y_k; \\ 2. \; i_{k+1} \in [1,...,n] \text{ – независимо от предыдущих розыгрышей;} \\ 3. \; y_{k+1} = \text{Grad}_{i_{k+1}}(x_{k+1}); \\ 4. \; z_{k+1} = \text{Mirr}_{z_k}\left(\alpha_{k+1} n \nabla_{i_{k+1}} f(x_{k+1})\right). \end{array}}$$

Оценка скорости сходимости такого метода

$$\boxed{N(\varepsilon) = \text{O}\left(n\sqrt{\frac{\Theta}{\varepsilon}} \ln\left(\frac{1}{\sigma}\right)\right),}$$

где

$$\boxed{\Theta = V_{x_0}(x_*).}$$

Получается эта оценка из следующей формулы (см. последнюю формулу в доказательстве леммы 4.3 [6])

$$\alpha_{k+1}^2 n^2 E_{i_{k+1}}\left[f(y_{k+1})\big|i_1,...,i_k\right] - \left(\alpha_{k+1}n^2 - \alpha_{k+1}\right)f(y_k) \le$$

$$\le \alpha_{k+1}\left\{f(x_{k+1}) + \left\langle \nabla f(x_{k+1}), u - x_{k+1}\right\rangle\right\} + V_{z_k}(u) - E_{i_{k+1}}\left[V_{z_{k+1}}(u)\big|i_1,...,i_k\right].$$

Взяв математическое ожидание $E_{i_1,...,i_k}[\cdot]$ (если $k \ge 1$) от каждого такого неравенства, и просуммировав то что получается по $k = 0,..., N-1$, получим

$$\alpha_N^2 n^2 E_{y_N}\left[f(y_N)\right] \le \min_{u \in Q}\left\{\sum_{k=0}^{N-1}\alpha_{k+1}E_{x_{k+1}}\left[f(x_{k+1}) + \left\langle \nabla f(x_{k+1}), u - x_{k+1}\right\rangle\right] + V_{z_0}(u) - E_{z_N}\left[V_{z_N}(u)\right]\right\} \le$$

$$\le \left(\sum_{k=0}^{N-1}\alpha_{k+1}\right)f_* + V_{z_0}(x_*) - E_{z_N}\left[V_{z_N}(x_*)\right].$$

Из последнего неравенства получается нужная оценка $N(\varepsilon)$ (только для сходимости в среднем, для получения оценки вероятностей больших уклонений необходимо использовать[10] неравенство концентрации Азума–Хефдинга для последовательностей мартингал-разностей, см., например, главу 7 [7]). Приведенное неравенство сразу показывает прямо-двойственность метода [31], [33], что это означает (и какая от этого польза), хорошо можно продемонстрировать конкретными примерами [5], [32], [35], [36] (см. также пример 3 ниже). Также из приведенной оценки следует, что

$$E_{z_k}\left[\frac{1}{2}\|z_k - x_*\|^2\right] = E_{z_k}\left[V_{z_k}(x_*)\right] \le E_{z_{k-1}}\left[V_{z_{k-1}}(x_*)\right] \le ... \le E_{z_{k-1}}\left[V_{z_1}(x_*)\right] \le V_{z_0}(x_*).$$

Можно привести и более точные вероятностные оценки на субмартингал

$$\|z_k - x_*\|^2.$$

Аналогичными субмартингальными свойствами обладают и последовательности[11]

$$\|y_k - x_*\|^2, \|x_k - x_*\|^2,$$

что доказывается по индукции, исходя из пп. 1, 3 в определении ACRCD*, выпуклости квадрата нормы и неравенства Йенсена. Детали мы вынуждены здесь опустить. В качестве "сухого остатка" можно сформулировать следующий результат (см. также [9], [32]): с вероятностью $\ge 1 - \sigma$

---

[10] Впрочем, можно получить неравенства на вероятности больших уклонений с помощью неравенства Маркова подобно тому, как это было описано выше в п. 4 для ACRCD.

[11] При доказательстве этого факта существенно используется евклидовость нормы, к сожалению, для не евклидовых норм похоже, что результат перестает быть верным (см., например, приложение B.1 в [6]) если по-прежнему исходить из ПГМ в представлении БГМ (не пытаясь заменить его, например, на ППГМ).

$$\max_{k=1,\ldots,N}\left\{\|y_k - x_*\|^2, \|z_k - x_*\|^2, \|x_k - x_*\|^2\right\} \le R^2,$$

где

$$R^2 = CV_{z_0}(x_*)\ln(N/\sigma),$$

а $C < 100$ – некоторая универсальная константа. По-сути, это означает, что если заранее знать $V_{z_0}(x_*)$, то, например, константы Липшица можно определять не на всем $Q$ (если $Q$ не ограничено, то и константы могут быть не ограничены), а на пересечении $Q$ с "шаром" в $\|\ \|$-норме с центром в точке $x_0$ и радиуса $R$. Вместе с прямо-двойственностью, это оказывается полезным инструментом для использования покомпонентных методов при решении двойственных задач [5] (см. также пример 3 ниже). Из описания ACRCD* также следует, что метод позволяет адаптивно подбирать константы Липшица по разным направлениям, подобно п. 6.1 работы [2].[12] Теперь уже нет проблемы с завышенной оценкой параметра $\Theta$, входящего в оценку $N(\varepsilon)$, в виду $\Theta = V_{x_0}(x_*)$. И из двух потенциально неизвестных априорно параметров $\Theta$ и $d$ теперь остается только один $\Theta$.

Можно показать, что в выписанной формуле для $N(\varepsilon)$ константа в $\mathrm{O}(\ )$ не больше, чем в теореме 2. Используя эту явную формул для $N(\varepsilon)$, подобно п. 5 работы [6], с помощью техники рестартов (по расстоянию от текущей точки до решения) можно перенести полученные результаты на случай $\mu$-сильно выпуклой в норме $\|\ \|$ функции (заметим, что при таком перенесении можно сохранить возможность метода адаптивно настраиваться на константы Липшица). Соответствующая оценка числа итераций будет иметь следующий вид (для евклидовой нормы $\|\ \|$, для неевклидовой под корнем может возникнуть дополнительный логарифмический по $n$ множитель [9])

$$N(\varepsilon) = \mathrm{O}\left(n\sqrt{\frac{1}{\mu}\ln\left(\frac{\ln(\mu\Theta/\varepsilon)}{\sigma}\right)}\ln\left(\frac{\mu\Theta}{\varepsilon}\right)\right).$$

**Замечание 3. (обобщение на блочно-компонентные методы и на более общие прокс-структуры и множества)** Описанный метод допускает следующее обобщение.[13] Пусть (см. также теорему[14] 5 [2])

---

[12] Чтобы сохранить дешевизну итерации (эффективность метода) для композитных постановок или в случае когда множество $Q$ не параллелепипедного типа здесь требуются некоторые оговорки (подобные сделанным в замечаниях 7, 8 ниже) о возможности эффективно пересчитывать значения функции.

[13] Мы не будем подробно пояснять все используемые далее обозначения – они стандартны и должны быть понятны из контекста, детали см., например, в [2], [3].

$$x = (x_1, \ldots, x_n), \quad Q = \prod_{i=1}^{n} Q_i.$$

Каждый $x_i \in Q_i$, в свою очередь, является вектором (размерности у этих векторов могут быть разными). Пусть в соответствующих подпространствах (отвечающих различным блокам) введены нормы $\left\{\sqrt{L_i} \|x_i\|_i\right\}_{i=1}^{n}$ и соответствующие этим нормам "расстояния" Брегмана $\left\{V_{x_i}^i(y_i)\right\}_{i=1}^{n}$ (см., например, [7], [9]). Положим

$$\|x\|^2 = \sum_{i=1}^{n} L_i \|x_i\|_i^2, \quad \|\nabla f(x)\|_*^2 = \sum_{i=1}^{n} L_i^{-1} \|\text{grad}_{x_i} f(x)\|_{i,*}^2, \quad V_x(y) = \sum_{i=1}^{n} V_{x_i}^i(y_i).$$

Будем считать, что для всех $x, x + h\tilde{e}_i \in Q$

$$\|\text{grad}_{x_i} f(x + h\tilde{e}_i) - \text{grad}_{x_i} f(x)\|_{*,i} \leq L_i h \|[\tilde{e}_i]_i\|_i,$$

где вектор $\tilde{e}_i$ имеет все нули в компонентах, не соответствующих $i$-му блоку. Введённые обозначения позволяют переписать сам метод (см. также [6]). При этом оценки будут иметь точно такой же вид, меняется только интерпретация параметров, норм (расстояний) в этих оценках. Заметим, что в приложениях к поиску равновесий в популяционных играх загрузок с большим числом популяций (в частности, задачах поиска равновесного распределения потоков по путям в графе транспортной сети [32], [35], [36]), часто возникают множества имеющие вид прямого произведения симплексов. До настоящего момента было не понятно, можно ли (а если можно, то как) применять к таким задачам покомпонентные методы.

**Замечание 4. (обобщение на задачи стохастической оптимизации)** Предположим, что исходная задача имеет вид

$$E_\xi [f(x;\xi)] \to \min_{x \in Q}.$$

Если $f(x;\xi)$ – выпуклая по $x$ функция (при всех $\xi$) с константами Липшица (равномерно не только по $x, x + he_i \in Q$, но и по $\xi$)

$$\|\text{grad}_{x_i} f(x + h\tilde{e}_i; \xi) - \text{grad}_{x_i} f(x; \xi)\|_{*,i} \leq L_i h \|[\tilde{e}_i]_i\|_i,$$

где вектор $e_i$ имеет все нули в компонентах, не соответствующих $i$-му блоку. Введём

---

[14] Отметим также, что везде в этой теореме можно вместо $R_1^2(x_0)$ писать $2\|x_0 - x_*\|_1^2$, что немного улучшает оценку теоремы.

$$D = \max_{x \in Q} E_\xi \left[ \left\| \nabla f(x;\xi) - E_\xi \left[ \nabla f(x;\xi) \right] \right\|_*^2 \right].$$

Тогда если вместо $\nabla_i f(x;\xi)$ можно рассчитывать только на $\nabla_i f(x;\xi)$:

$$E_\xi \left[ \nabla_i f(x;\xi) \right] \equiv \nabla_i E_\xi \left[ f(x;\xi) \right],$$

то оценка в теореме 2 изменится следующим образом (см. также [9])

$$N(\varepsilon) := O\left( \max\left\{ N(\varepsilon), n \frac{D\Theta}{\varepsilon^2} \ln\left( \frac{1}{\sigma} \right) \right\} \right),$$

в случае $\mu$-сильно выпуклой в норме $\|\ \|$ функции

$$N(\varepsilon) := O\left( \max\left\{ N(\varepsilon), n \frac{D}{\mu\varepsilon} \ln\left( \frac{\ln(N(\varepsilon))}{\sigma} \right) \right\} \right).$$

**Замечание 5. (учет ошибок в вычислении компонент градиента)** Немного специфицируя концепцию неточного оракула ($\delta \geq 0$ – уровень шума) из главы 4 [7], введем следующее предположение (векторы $\tilde{e}_i$, $\nabla_i f_{\delta,L_i}(x)$ имеет все нули в компонентах, не соответствующих $i$-му блоку): для любого $x \in Q$ существуют такие $f_{\delta,L_i}(x)$ и $\nabla_i f_{\delta,L_i}(x)$, что для всех $y = x + h\tilde{e}_i \in Q$ выполняется

$$0 \leq f(y) - f_{\delta,L_i}(x) + \langle \nabla_i f_{\delta,L_i}(x), y - x \rangle \leq \frac{L_i}{2} \left\| [y-x]_i \right\|_i^2 + \delta,$$

т.е.

$$0 \leq f(x + h\tilde{e}_i) - f_{\delta,L_i}(x) + h\langle \nabla_i f_{\delta,L_i}(x), \tilde{e}_i \rangle \leq \frac{L_i h^2}{2} \left\| [\tilde{e}_i]_i \right\|_i^2 + \delta.$$

При $\delta = 0$ отсюда получаем определение констант Липшица $L_i$ в блочно-покомпонентном методе из замечания 3. Полезными следствиями введенного определения являются следующие неравенства

$$0 \leq f(x) - f_{\delta,L_i}(x) \leq \delta,$$

$$\frac{1}{2} \left\| \nabla_i f_{\delta,L_i}(x) \right\|_*^2 \leq f(x) - f\left( \text{Grad}_i^\delta(x) \right) + \delta,$$

где

$$\text{Grad}_i^\delta(x) = \arg\min_{\tilde{x} \in Q} \left\{ \langle \nabla_i f_{\delta,L_i}(x), \tilde{x} - x \rangle + \frac{1}{2} \|\tilde{x} - x\|^2 \right\} = x - \frac{1}{L_i} \nabla_i f_{\delta,L_i}(x),$$

здесь для наглядности мы ограничились случаем евклидовой нормы и $Q = \mathbb{R}^n$. С помощью этих неравенств можно скорректировать (см. также [9]) оценку теоремы 2 (и различные ее обобщения) на случай, когда мы можем вычислять вместо "честных" компонент градиента $\nabla_i f(x)$ только их приближенные (в указанном выше смысле) аналоги $\nabla_i f_{\delta, L_i}(x)$:

$$\boxed{f\left(x^{N(\varepsilon)}\right) - f_* \leq \varepsilon + \mathrm{O}(N(\varepsilon)\delta).}$$

Вообще эта формула типична для всех известных нам ускоренных методов (полноградиентных, покомпонентных, прямых). И это соответствует самому худшему (быстрому) варианту накопления ошибки. Для неускоренных методов $\mathrm{O}(N\delta) \to \mathrm{O}(\delta)$, что соответствует самому лучшему варианту накопления ошибки (то есть когда такое накопление отсутствует). Подобно полноградиентным методам для покомпонентных методов (и прямых) можно предложить так называемы промежуточные методы (см., например, [7]–[9]) с накоплением ошибки $\mathrm{O}(N^p \delta)$, $p \in [0,1]$. Из этого замечания (а также замечания 2) возникает гипотеза о возможности создания универсального покомпонентного метода [28].

Интересно было бы объединить замечания 4, 5 с целью получения покомпонентной версии результатов главы 7 [7] и [8].

**Замечание 6. (обобщение на взвешенную рандомизацию)** Предположим, что вместо одной компоненты можно выбирать блок компонент (всего $n$ блоков), причем, вообще говоря, с разными вероятностями: выбираем блок компонент $i$ с вероятностью[15]

$$\boxed{p_i = \frac{L_i^\beta}{\sum_{j=1}^n L_j^\beta}, \ i = 1, \ldots, n,}$$

где параметр степени $\beta \in [0,1]$. При этом необходимо будет переопределить норму

$$\boxed{\|x\|^2 = \sum_{i=1}^n L_i x_i^2 \to \sum_{i=1}^n L_i^{1-2\beta} \|x_i\|_i^2,}$$

а, соответственно, также прокс-функцию и параметр $\Theta$. При этом во всех приведенных выше формулах, которые определяют метод, в частности, для ACRCD это

$$\tau = \frac{1}{\alpha n^2 + 1}, \ \alpha = \frac{1}{n}\sqrt{\frac{\Theta}{d}}, \ K = 8n\sqrt{\frac{\Theta}{d}}, \ N = 27n\sqrt{\frac{\Theta}{\varepsilon}} \log_2\left(\frac{\log_2(d/\varepsilon)}{\sigma}\right),$$

необходимо будет сделать замену

---

[15] Приготовление памяти для генерирования из описанного распределения стоит $\mathrm{O}(n)$. Это делается один раз (строится соответствующее двоичное дерево Л.В. Канторовича [2]). Случайные разыгрывания $i$ при наличии правильно подготовленной памяти будут стоить $\mathrm{O}(\log_2 n)$ – каждое.

$$\boxed{n \to \sum_{i=1}^{n} L_i^{\beta}.}$$

Выше в этом пункте мы рассматривали случай $\beta = 0$ (в другом ключе этот случай также рассматривался в [2]). Можно ожидать, что на практике этот вариант предпочтительнее. Заметим, что ранее уже рассматривались отдельно случаи $\beta = 1/2$ [27] и $\beta = 1$ [26].

В связи с замечанием 6 возникает вопрос: существуют ли еще более общие способы (с большим числом степеней свободы) сочетания выбора рандомизации и нормы? Положительный ответ, более менее, очевиден (см. также [4]), но интересно было бы предложить такие способы, которые в определенных ситуациях позволяли бы еще более ускориться по сравнению с методом, порожденным замечаниями 2, 6.

**Замечание 7. (стоимость итерации / не разреженный случай)** Рассмотрим, следуя Ю.Е. Нестерову (см. слайд 7 [27]), следующий случай

$$f(x) = F(Ax, x), \; x \in \mathbb{R}^n, \; y = Ax \in \mathbb{R}^m.$$

Будем считать, что значение $F(y, x)$ (а, следовательно, и градиент [38]) можно посчитать за $\mathrm{O}(m+n)$. Пусть верно хотя бы одно из следующих условий: 1) $n = \mathrm{O}(m)$; 2) расчет $\mathrm{grad}_y F(y,x)$ стоит $\mathrm{O}(m)$, а $\partial F(y,x)/\partial x_j$ – $\mathrm{O}(m)$.[16] Тогда амортизационная (средняя) сложность одной итерации будет $\mathrm{O}(m)$. Обоснование этого факта можно получить как простое следствие более общих рассуждений, проводимых в следующем замечании.

**Замечание 8. (стоимость итерации / разреженный случай)** Из первого пункта описания алгоритма ACRCD кажется, что всегда один шаг этого алгоритма будет требовать, как минимум $\geq n$ арифметических операций. Однако замечание 7 показывает (при $m \ll n$), что это совсем не обязательно. Естественно возникает вопрос: а можно ли получить еще больше (еще более дешевую итерацию)? В определенных (разреженных задачах специальной структуры) ответ оказывается положительным. Пояснению этого тезиса и будет посвящена оставшаяся часть данного замечания. Оказывается, что при наличии у задачи определенной структуры (например, в случае $f(x) = x^T A x$ или $f(x) = \|Ax - b\|_2^2$), нет необходимости выполнять первый пункт честно (в полном объеме). Далее мы описываем идею, заимствованную из работ [3], [26], [27]. Предварительно перепишем алгоритм

---

[16] В ряде приложений посчитать одну компоненту градиента оказывается в $n$ раз дешевле, чем сам градиент. Например, это так для $f(x) = x^T A x$. Но все же верно это далеко не всегда. Например, для функции (см. пример 3 п. 5 ниже)

$$f(x) = \ln\left(\sum_{k=1}^{n} \exp(x_k)\right)$$

стоимость расчета самой функции, ее градиента и любой компоненты градиента одинаковы по порядку.

ACRCD в рассматриваемом нами случае (рассматривается задача безусловной оптимизации) следующим (эквивалентным) образом ($x_0 = u_0 = v_0$)

$$\text{ACRCD}'(\alpha, \tau; \Theta, x_0, f(x_0) - f_* \leq d)$$

1. $x_{k+1} = (1-\tau)^{k+1} v_k + u_k$;

2. $i_{k+1} \in [1,...,n]$ – независимо от предыдущих розыгрышей;

3. $v_{k+1} = v_k + \dfrac{1}{L_{i_{k+1}}} \dfrac{\alpha n - 1}{(1-\tau)^k} \nabla_{i_{k+1}} f(x_{k+1})$;

4. $u_{k+1} = u_k - \dfrac{\alpha n}{L_{i_{k+1}}} \nabla_{i_{k+1}} f(x_{k+1})$.

5. если $k = N-1$, то выдаем

$y_N = y_{k+1} = \text{Grad}_{i_{k+1}}(x_{k+1}) = \text{Grad}_{i_N}(x_N)$.

Предположим, что

$$f(x) = \sum_{r=1}^{m} \varphi_r(a_r^T x),$$

где функции $\varphi_r$ – простой структуры, т.е. дифференцируемые за $\text{O}(1)$ каждая, причем $A = \|a_1 \dots a_m\|^T$ – разреженная матрица (число ненулевых элементов $sn$, т.е. в каждом столбце в среднем $s \ll m$ ненулевых элементов). Тогда подобно [3] один шаг метода ACRCD' (кроме самого первого шага, который может стоить $\text{O}(n)$) может быть осуществлен в среднем за $\text{O}(s)$ (амортизационная сложность). Действительно, пункт 2 ACRCD' можно осуществить за $\text{O}(\ln n)$ (считаем $\ln n = \text{O}(s)$). Пункты 3, 4 за $\text{O}(s)$. Такая стоимость этих пунктов обусловлена необходимостью пересчета $\partial f(x_{k+1})/\partial x_i$. Осуществлять этот пересчет необходимо, используя пункт 1. Покажем, как можно это эффективно делать. Прежде всего, заметим, что если мы уже посчитали $a_r^T v_k$ и $a_r^T u_k$, то посчитать дополнительно $a_r^T x_{k+1}$ будет стоить $\text{O}(1)$. Также заметим, что если мы уже посчитали $Av_k$ и $Au_k$, то посчитать дополнительно $Av_{k+1}$ и $Au_{k+1}$ будет стоить $\text{O}(s)$. Чтобы посчитать $\partial f(x_{k+1})/\partial x_i$ нужно вычислить частные производные (по $x_i$) в среднем у $s$ слагаемых в сумме. Каждая такая частная производная (по предположению о простоте структуры функций $\varphi_r$) рассчитывается за $\text{O}(1)$, в предположении известности всех аргументов этих функций. Учитывая, что на пересчет всех аргументов уйдет $\text{O}(s)$ (см. выше), то в сред-

нем, общие трудозатраты будут $\mathrm{O}(s)+\mathrm{O}(s)=\mathrm{O}(s)$. Таким образом, один шаг метода будет стоить $\mathrm{O}(s)$. Чтобы посчитать выход алгоритма:

$$\bar{x}_K = \frac{1}{K}\sum_{k=1}^{K} x_k = \frac{1}{K}\sum_{k=1}^{K}\sum_{\tilde{K}=k}^{K}\left(\frac{\alpha n-1}{L_{i_k}}\frac{(1-\tau)^{\tilde{K}}}{(1-\tau)^k} - \frac{\alpha n}{L_{i_k}}\right)\nabla_{i_k} f(x_k) =$$

$$= \sum_{k=1}^{K}\frac{1}{K}\left(\frac{\alpha n-1}{L_{i_k}}\frac{1-\tau}{\tau}\left(1-(1-\tau)^{K-k}\right) - \frac{\alpha n}{L_{i_k}}(K-k+1)\right)\nabla_{i_k} f(x_k),$$

достаточно $\mathrm{O}(K)$ арифметических операций. Все написанное выше переносится и на ACRCD*.

$$\text{ACRCD*}'(x_0 = u_0 = v_0 = w_0)$$

1. $x_{k+1} = \left\{\prod_{j=1}^{k}(1-\tau_j)\right\}\cdot\left[u_k + v_k + \left\{\sum_{j=1}^{k}\frac{\tau_j}{\prod_{l=1}^{j}(1-\tau_l)}\right\}w_k\right];$

2. $i_{k+1} \in [1,...,n]$ – независимо от предыдущих розыгрышей;

3. $u_{k+1} = u_k + \frac{\alpha_{k+1}n}{L_{i_{k+1}}}\sum_{j=1}^{k}\frac{\tau_j}{\prod_{l=1}^{j}(1-\tau_l)}\nabla_{i_{k+1}} f(x_{k+1});$

4. $v_{k+1} = v_k - \frac{1}{L_{i_{k+1}}}\frac{1}{\prod_{j=1}^{k}(1-\tau_j)}\nabla_{i_{k+1}} f(x_{k+1});$

5. $w_{k+1} = w_k - \frac{\alpha_{k+1}n}{L_{i_{k+1}}}\nabla_{i_{k+1}} f(x_{k+1});$

6. если $k = N-1$, то выдаем
$y_N = y_{k+1} = \mathrm{Grad}_{i_{k+1}}(x_{k+1}) = \mathrm{Grad}_{i_N}(x_N).$

## 5. Примеры применения ускоренных покомпонентных методов

Начнем этот пункт с очень простого примера, демонстрирующего, что покомпонентные методы (в том числе неускоренные), вообще говоря, не применимы к произвольной задаче выпуклой оптимизации.

**Пример 1.** Рассмотрим выпуклую задачу

$$f(x) = (x_1 - 2)^2 + (x_2 - 1)^2 \to \min_{x \in Q},$$

$$Q = \{x = (x_1, x_2) \geq 0 : \quad x_1 + x_2 \leq 2\}.$$

Предположим, что покомпонентный метод стартует с точки

$$x^0 = (1, 1) \in Q.$$

Тогда метод (если в методе жестко прописано оптимизировать по выбранному направлению – у нас это не так) не сдвинется с места по обоим направлениям, задаваемым ортами и проходящим через $x^0$ (поскольку $f(x)$ внутри $Q$ имеет минимум в $x^0$ по этим направлениям), в то время как в точке

$$x^* = (0.5, 1.5) \in Q$$

$f(x)$ достигает минимума на $Q$:

$$f(x^*) = 0.5 < 1 = f(x^0).$$

Если в методе жестко прописано оптимизировать по выбранному направлению (в подавляющем большинстве существующих вариантов покомпонентных методов именно так и сделано), то для возможности рассмотрения случая $Q \neq \mathbb{R}^n$ нужно дополнительно потребовать, чтобы для любого $i = 1, ..., n$ имело место условие: для любого $x = (x_1, ..., x_n) \in Q$ выполняется

$$(x_1, ..., x_{i-1}, x_{*i}, x_{i+1}, ..., x_n) \in Q.$$

Если множество $Q = \prod_{k=1}^{n}[a_k, b_k]$ – сепарабельно, то это условие, очевидно, выполняется. □

Перейдем теперь непосредственно к приложениям покомпонентных методов. За более подробной информацией о покомпонентных методах и примерах их приложений можно рекомендовать обратиться к [17], [25], [39].

**Пример 2.** Возьмем функцию (в изложении этого примера мы во многом следуем Ю.Е. Нестерову, см. слайд 8 [27])

$$f(x) = \frac{1}{2}\langle x, Sx \rangle - \langle b, x \rangle,$$

где $S$ – симметричная матрица, все элементы которой числа от 1 до 2. Возьмем метод ACRCD* в варианте замечания 6 с $\beta = 1/2$. Выберем евклидову норму. Константа Липшица этой функции по определению есть

$$L = \lambda_{\max}(S) \geq \lambda_{\max}\left(1_n 1_n^T\right) = n.$$

При этом покомпонентный метод дает константы Липшица $L_i = S_{ii} \leq 2$. Таким образом, получаем ускорение в $\sim \sqrt{n}$ раз. Действительно (см. замечание 6),

$$n \to \sum_{i=1}^{n} \sqrt{L_i} \leq \sqrt{2} n,$$

поэтому оценка числа итераций соответствующего покомпонентного быстрого градиентного метода (ПБГМ) будет

$$N_{ПБГМ}(\varepsilon) = O\left(n\sqrt{\frac{\Theta}{\varepsilon}} \ln\left(\frac{1}{\sigma}\right)\right),$$

а стоимость одной итерации $O(n)$. Итого

$$T_{ПБГМ} = O\left(n^2 \sqrt{\frac{\Theta}{\varepsilon}} \ln\left(\frac{1}{\sigma}\right)\right).$$

Для обычного (не покомпонентного) БГМ соответствующая оценка числа итераций имеет вид (отметим, что при выборе $\beta = 1/2$ можно считать $\Theta$ в обеих формулах одинаковым)

$$N_{БГМ}(\varepsilon) = O\left(\sqrt{\frac{L\Theta}{\varepsilon}}\right) = O\left(\sqrt{\frac{n\Theta}{\varepsilon}}\right).$$

Зато одна итерация стоит $O(n^2)$. Итого

$$T_{БГМ} = O\left(n^{5/2} \sqrt{\frac{\Theta}{\varepsilon}}\right).$$

В более общем случае полезно иметь в виду следующие неравенства

$$\frac{1}{n}\operatorname{tr}(S) \leq \lambda_{\max}(S) \leq \operatorname{tr}(S), \quad \frac{1}{n}\sum_{i=1}^{n}\sqrt{L_i} \leq \sqrt{\frac{1}{n}\sum_{i=1}^{n} L_i} = \sqrt{\frac{1}{n}\operatorname{tr}(S)}.$$

Таким образом (здесь мы опустили логарифмический множитель в оценке $T_{БГМ}$, поэтому вместо $O(\ )$ ввели $\tilde{O}(\ )$),

$$T_{\text{ПБГМ}} = \tilde{O}\left(n^2 \sqrt{\frac{(\operatorname{tr}(S)/n)\Theta}{\varepsilon}}\right) \le O\left(n^2 \sqrt{\frac{\lambda_{\max}(S)\Theta}{\varepsilon}}\right) = T_{\text{БГМ}}.$$

В разреженном случае, согласно замечанию 8, пропорции сохраняются

$$T_{\text{ПБГМ}} = \tilde{O}\left(sn \sqrt{\frac{(\operatorname{tr}(S)/n)\Theta}{\varepsilon}}\right) \le O\left(sn \sqrt{\frac{\lambda_{\max}(S)\Theta}{\varepsilon}}\right) = T_{\text{БГМ}}.$$

Обратим внимание, что выгода в $\sim \sqrt{n}$ раз является максимально возможной. Достигается она в ситуациях, когда $\lambda_{\max}(S)$ и $\operatorname{tr}(S)$ одного порядка. Скажем, если собственные значения матрицы $S$: $\{1,...,n\}$, то $\lambda_{\max}(S) = n$, а $\operatorname{tr}(S) \sim n^2$, т.е. нужна большая (более резкая) асимметрия. Если под матрицей $S$ понимать гессиан функционала задачи в "худшей" (с точки зрения рассматриваемых оценок) точке, то выписанные формулы не изменятся. Однако в разреженном случае потребуются большие оговорки, чтобы можно было сполна учесть разреженность в стоимости итерации. Поскольку в приложениях довольно типично выполнение неравенства

$$\frac{1}{n}\operatorname{tr}(S) \ll \lambda_{\max}(S),$$

то из приведенных оценок следует, что во многих случаях, получается ускорить вычисления за счет использования ПБГМ вместо БГМ (не говоря уже о возможности распараллеливания [3], [4]). Как уже отмечалось, в ряде случаев это выгода может достигать $\sim \sqrt{n}$ раз. Другие примеры, когда похожие пропорции имеют место можно посмотреть в работах [5], [32]. □

**Пример 3.** Рассмотрим следующую задачу энтропийно-линейного программирования (см., например, слайд 9 [27] и [37])

$$f(x) = \sum_{i=1}^{n} x_i \ln(x_i) \to \min_{x \in S_n(1);\, Ax=b},$$

$$S_n(1) = \left\{ x \in \mathbb{R}^n : x_i \ge 0,\ i=1,...,n,\ \sum_{i=1}^{n} x_i = 1 \right\},$$

причем будем считать (в связи с различными транспортными приложениями это представляется довольно естественным [32], [35]–[37]), что условие $\sum_{i=1}^{n} x_i = 1$ является следствием системы $Ax = b$. Построим двойственную задачу

$$\min_{x \in S_n(1);\, Ax=b} \sum_{i=1}^{n} x_i \ln(x_i) = \min_{x \in S_n(1)} \max_{y \in \mathbb{R}^m} \left\{ \sum_{i=1}^{n} x_i \ln(x_i) + \langle y, b - Ax \rangle \right\} =$$

$$= \max_{y \in \mathbb{R}^m} \min_{x \in S_n(1)} \left\{ \sum_{i=1}^{n} x_i \ln(x_i) + \langle y, b - Ax \rangle \right\} = \max_{y \in \mathbb{R}^m} \left\{ \langle y, b \rangle - \ln\left( \sum_{i=1}^{n} \exp\left( \left[ A^T y \right]_i \right) \right) \right\}.$$

Но с учетом написанного выше, двойственную задачу можно строить и по-другому

$$\min_{x \in S_n(1);\, Ax=b} \sum_{i=1}^{n} x_i \ln(x_i) = \min_{Ax=b} \sum_{i=1}^{n} x_i \ln(x_i) = \min_{x \in \mathbb{R}_+^n} \max_{y \in \mathbb{R}^m} \left\{ \sum_{i=1}^{n} x_i \ln(x_i) + \langle y, b - Ax \rangle \right\} =$$

$$= \max_{y \in \mathbb{R}^m} \min_{x \in \mathbb{R}_+^n} \left\{ \sum_{i=1}^{n} x_i \ln(x_i) + \langle y, b - Ax \rangle \right\} = \max_{y \in \mathbb{R}^m} \left\{ \langle y, b \rangle - \sum_{i=1}^{n} \exp\left( \left[ A^T y \right]_i - 1 \right) \right\}.$$

В обоих случаях решение прямой задачи можно восстановить по решению двойственной:

1) $\displaystyle x_i(y) = \frac{\exp\left( \left[ A^T y \right]_i \right)}{\sum_{k=1}^{n} \exp\left( \left[ A^T y \right]_k \right)}$, $i = 1, \dots, n$,

2) $x_i(y) = \exp\left( \left[ A^T y \right]_i - 1 \right)$, $i = 1, \dots, n$.

Таким образом, можно работать с двумя различными двойственными задачами:

1) $\displaystyle \varphi_1(y) = \ln\left( \sum_{i=1}^{n} \exp\left( \left[ A^T y \right]_i \right) \right) - \langle y, b \rangle \to \min_{y \in \mathbb{R}^m}$,

2) $\displaystyle \varphi_2(y) = \sum_{i=1}^{n} \exp\left( \left[ A^T y \right]_i - 1 \right) - \langle y, b \rangle \to \min_{y \in \mathbb{R}^m}$.

Исходя из явного двойственного (Лежандрова) представления этих функционалов, можно вычислить константу Липшица градиента и соответствующие константы Липшица по направлениям. Во втором случае, к сожалению, константы получаются не ограниченными, а вот в случае 1 они ограничены,[17] и их можно оценить, соответственно, как:

$$L_{БГМ} = \max_{k=1,\dots,n} \left\| A^{\langle k \rangle} \right\|_2^2 \quad \text{и} \quad L_{ПБГМ}^k \le L_{ПБГМ} = \max_{\substack{i=1,\dots,m \\ j=1,\dots,n}} \left| A_{ij} \right|^2 \le L_{БГМ},\ k=1,\dots,n.$$

---

[17] За счет сильной выпуклости энтропии в 1-норме с константой 1 на единичном симплексе (на положительном ортанте энтропия строго выпукла, но не сильно выпукла), из теоремы 1 [40] имеем:

$$\left\| \nabla \varphi_1(y_2) - \nabla \varphi_1(y_1) \right\|_q = \left\| Ax(y_2) - Ax(y_1) \right\|_q \le \left\| A^T \right\|_{p,1}^2 \left\| y_2 - y_1 \right\|_p,\ \left\| A^T \right\|_{p,1}^2 = \max_{\|y\|_p \le 1,\, \|x\|_1 \le 1} \langle A^T y, x \rangle^2 = \max_{k=1,\dots,n} \left\| A^{\langle k \rangle} \right\|_q^2,$$

где $1/p + 1/q = 1$. Беря $p = 2$, получим константу Липшица градиента $\varphi_1(y)$: $\displaystyle \max_{k=1,\dots,n} \left\| A^{\langle k \rangle} \right\|_2^2$. Беря $p = 1$ можно получить, что константа Липшица производной $\varphi_1(y)$ по каждому направлению не больше $\displaystyle \max_{\substack{i=1,\dots,m \\ j=1,\dots,n}} \left| A_{ij} \right|^2$.

Будем считать, что все элементы матрицы $A_{ij}$ удовлетворяют условию: $1 \le A_{ij} \le 2$. Тогда $L_{БГМ} \ge m$, а $L_{ПБГМ} \le 4$. Выберем в двойственном пространстве евклидову прокс-структуру.[18] Решая первую двойственную задачу БГМ (стартуем в точке 0), получим следующую оценку времени работы метода

$$T^1_{БГМ} = O\left(mn\sqrt{\frac{L_{БГМ}\Theta}{\varepsilon}}\right) = O\left(mn\sqrt{\frac{m\Theta}{\varepsilon}}\right).$$

Если же применить к первой двойственной задаче ACRCD* (стартуем в точке 0) с $\beta = 1/2$ или $\beta = 0$ (см. замечания 2, 6, 7), то получим, что с вероятностью $\ge 1 - \sigma$

$$T^1_{ПБГМ} = \tilde{O}\left(nm\sqrt{\frac{L_{ПБГМ}\Theta}{\varepsilon}}\right) = \tilde{O}\left(mn\sqrt{\frac{\Theta}{\varepsilon}}\right).$$

В обеих формулах $\Theta$ – квадрат евклидового размера двойственного решения. Таким образом, за счет использования ПБГМ удается ускориться в $\sim \sqrt{m}$ раз. Все изложенное в этом примере распространяется и на случай когда вместо ограничений в виде равенств (или наряду с ними) мы имеем ограничения в виде неравенств $Cx \le d$. Если по-прежнему обозначать общее число ограничений через $m$, то выписанные формулы останутся справедливыми. Пока мы решили только двойственную задачу с точностью по функции $\varepsilon$. То есть в двойственном пространстве используемый нами метод сгенерировал последовательности $\{y_k\}_{k=1}^N$, $\{z_k\}_{k=1}^N$, $\{\alpha_k\}_{k=1}^N$ (см. обозначения замечания 2 с заменой $x_k \to z_k$, сделанной во избежание путаницы) такие, что (для БГМ математическое ожидание можно не писать)

$$E_{y_N}\left[\varphi_1(y_N)\right] - \varphi_* \le \varepsilon, \; \varphi_* = f_*.$$

Оказывается, что не использующаяся в этой формуле и накопленная методом информация $\{z_k\}_{k=1}^N$, $\{\alpha_k\}_{k=1}^N$ позволяет восстанавливать с такой же точностью решение прямой задачи (детали см., например, в работе [5]):

$$\bar{x}^N = \frac{1}{S_N}\sum_{k=1}^N \alpha_k x(z_k), \; S_N = \sum_{k=1}^N \alpha_k,$$

$$\sqrt{\Theta}\left\|A\bar{x}^N - b\right\|_2 \le \varepsilon, \; \left|f(\bar{x}^N) - f_*\right| \le \varepsilon.$$

Несложно показать, что можно так организовать вычисление $\bar{x}^N$, что выписанные ранее оценки трудоемкости методов БГМ и ПБГМ (в категориях O( )) не изменятся. Если те-перь рассмотреть разреженный случай ($s \ll m$ – среднее число ненулевых элементов в столбце матрицы $A$, $\tilde{s} = sn/m$ – в строчке), то оценка БГМ улучшится

$$T_{БГМ}^1 = O\left( sn\sqrt{\frac{L_{БГМ}\Theta}{\varepsilon}} \right),$$

в то время как оценка ПБГМ останется неизменной (см. замечание 8)

$$T_{ПБГМ}^1 = \tilde{O}\left( mn\sqrt{\frac{L_{ПБГМ}\Theta}{\varepsilon}} \right).$$

---
[18] Вообще при решении двойственных задач это совершенно естественно [37], поскольку оптимизация происходит либо на всем пространстве, либо на прямом произведении пространства и неотрицательно ортанта.

Получилось это из-за наличия (при построении двойственной задачи) связывающего переменные (симплексного) ограничения – не имеющего полностью сепарабельную структуру (то есть не распадающегося в прямое произведение ограничений на отдельные компоненты). Другое дело, если мы рассмотрим вторую двойственную задачу. Она полностью сепарабельная (подходит под замечание 7 в этом смысле). То есть для второй двойственной задачи можно найти ее решение с точностью по функции $\varepsilon$:

$$E_{y_N}\left[ \varphi_2(y_N) \right] - \varphi_* \leq \varepsilon$$

за время (обратим внимание, что $\overline{L}_{ПБГМ}$, $\overline{\Theta}$ соответствуют функционалу и решению второй двойственной задачи, и, вообще говоря, отличаются от введенных ранее $L_{ПБГМ}$, $\Theta$)

$$T_{ПБГМ}^2 = \tilde{O}\left( \tilde{s}m\sqrt{\frac{\overline{L}_{ПБГМ}\overline{\Theta}}{\varepsilon}} \right) = \tilde{O}\left( sn\sqrt{\frac{\overline{L}_{ПБГМ}\overline{\Theta}}{\varepsilon}} \right).$$

При этом следует использовать метод ACRCD* из замечания 2 (если смотреть через замечание 6, то следует полагать $\beta = 0$) с адаптивным подбором констант Липшица, поскольку ограничить их не представляется возможным в виду априорного отсутствия информации о локализации решения двойственной задачи. Тем не менее, согласно замечанию 2 можно быть уверенным, что несмотря на неограниченность множества, на котором происходит оптимизация, и неограниченности констант Липшица на этом множестве, существуют "эффективные" константы, на которые мы адаптивно настраиваемся по ходу работы метода. К сожалению, восстановить за то же по порядку время решение прямой задачи $\bar{x}^N$ в данном случае не получается, даже если пытаться использовать соответствующие наработки замечания 8 по расчету $\bar{x}_K$. Мы снова возвращаемся к оценке типа

$$T_{ПБГМ}^2 = \tilde{O}\left( nm\sqrt{\frac{\overline{L}_{ПБГМ}\overline{\Theta}}{\varepsilon}} \right).$$

Если отказаться от ускоренности метода, то обычный неускоренный покомпонентный метод (ПМ) [1], [2] позволяет сохранить (с учетом необходимости восстановления решения прямой задачи) дешевую итерацию[19] $\mathrm{O}(\tilde{s})$

$$T_{\textit{ПМ}}^2 = \tilde{\mathrm{O}}\!\left( m\tilde{s}\,\frac{\overline{L}_{\textit{ПБГМ}}\overline{\Theta}}{\varepsilon} \right) = \tilde{\mathrm{O}}\!\left( sn\,\frac{\overline{L}_{\textit{ПБГМ}}\overline{\Theta}}{\varepsilon} \right).$$

В данном случае (и это довольно типично), этому методу скорее стоит предпочесть БГМ с итоговой оценкой

$$T_{\textit{БГМ}}^2 = \tilde{\mathrm{O}}\!\left( sn\sqrt{\frac{\overline{L}_{\textit{БГМ}}\overline{\Theta}}{\varepsilon}} \right),$$

но в таком случае уж лучше применять БГМ к первой двойственной задачи, обладающей лучшими свойствами. Получается, что необходимость восстановления решения прямой задачи для ПБГМ накладывает дополнительные ограничения на структуру задачи, чтобы можно было полноценно воспользоваться разреженностью. К сожалению, похоже, что эти дополнительные ограничения, фактически, оставляют возможность только для задач вида (и небольших "аффинных" релаксаций, например, добавление разреженных аффинных неравенств)

$$f(x) = \frac{1}{2}\|x - x_g\|_2^2 \to \min_{Ax=b},$$

полноценно использовать в разреженном случае описанный выше подход. При этом решение прямой и двойственной задачи можно восстанавливать по описанному выше механизму (с учетом линейности зависимости $x(y)$ удается воспользоваться техникой пересчета $\overline{x}_K$ из замечания 8 для восстановления $\overline{x}^N$). При этом оценки, соответствующих методов БГМ и ПБГМ, применённых к двойственной задаче (и стартующих с точки 0), примут следующий вид[20]

$$\tilde{T}_{\textit{БГМ}} = \tilde{\mathrm{O}}\!\left( sn\sqrt{\frac{\tilde{L}_{\textit{БГМ}}\tilde{\Theta}}{\varepsilon}} \right) \text{ и } \tilde{T}_{\textit{ПБГМ}} = \tilde{\mathrm{O}}\!\left( sn\sqrt{\frac{\tilde{L}_{\textit{ПБГМ}}\tilde{\Theta}}{\varepsilon}} \right).$$

---

[19] Это типично для неускоренных покомпонентных методов (в том числе прямо-двойственных), т.е. в отличие от ускоренных методов, тут требуется намного более слабые предположения, чтобы обеспечить выполнение условия: стоимость итерации покомпонентного метода дешевле стоимости итерации соответствующего полноградиентного метода в число раз по порядку равному размерности пространства.

[20] Используем следующие обозначения: $\tilde{\Theta}$ – квадрат евклидового размера решения двойственной задачи,

$$\tilde{L}_{\textit{БГМ}} = \max_{\|y\|_2\le 1, \|x\|_2\le 1}\left\langle A^T y, x\right\rangle^2 = \max_{\|x\|_2\le 1}\|Ax\|_2^2 = \lambda_{\max}\!\left(A^T A\right) \text{ и } \tilde{L}_{\textit{ПБГМ}} \le \max_{\|y\|_1\le 1, \|x\|_2\le 1}\left\langle A^T y, x\right\rangle^2 = \max_{\|x\|_2\le 1}\|Ax\|_\infty^2 = \max_{k=1,\ldots,m}\|A_k\|_2^2.$$

Написанное выше может навести на мысль, что ускоренные покомпонентные методы для двойственной задачи, как правило, не позволяют учитывать разреженность задачи. На самом деле это не так. За небольшую дополнительную плату (логарифмический множитель) можно специальным образом регуляризовать двойственную задачу (с помощью техники рестартов подобрать правильный параметр регуляризации, см. глава 3 [7], [37]), и использовать ПБГМ для регуляризованного функционала двойственной задачи, т.е. использовать ПБГМ в сильно выпуклом случае (см. текст сразу после замечания 2). При таком подходе достаточно просто решить (с желаемой точностью) двойственную задачу, а решение прямой задачи (в том же смысле, что и выше – с той же точностью) получается просто при подстановке найденного решения двойственной задачи в формулу $x(y)$. Тем не менее, здесь необходимо оговориться, что хотя описанный только что прием и "спасает положение", все же получается это за упомянутую дополнительную плату. Хотя по теоретическим оценкам это плата, действительно, небольшая, численные эксперименты показывают, что реальные потери при использовании такой регуляризации вместо прямодвойственности могут быть существенны. □

Результаты, изложенные в примере 3, допускают серьезные обобщения. В частности, можно переносить (частично), изложенное в примере 3, на сепарабельные задачи типа проектирования на аффинное многообразие

$$f(x) = \sum_{i=1}^{n} f_i(x_i) \to \min_{\substack{Ax=b,\, Cx \leq d \\ x \in Q_x}}$$

и более общий класс сепарабельных задач лежандровского типа (включающий проекционный класс)

$$\max_{x \in Q_x} \left\{ \langle y, Ax \rangle - \sum_{i=1}^{n} f_i(x_i) \right\} + \tilde{g}(y) \to \min_{y \in Q_y}.$$

Здесь $Q_x$ – множество простой структуры (в смысле проектирования), $Q_y$ – множество, подходящее для эффективного использования (блочно-)покомпонентных методов (см. п. 4), $\tilde{g}(y)$ – "хорошая" функция для покомпонентных методов (см. п. 4). Можно также не предполагать явной формулы, связывающей $x(y)$ (тогда потребуется еще воспользоваться замечанием 5). В определенных ситуациях можно даже пытаться отказаться от сепарабельности $f(x)$ (к сожалению, вот тут пока мало что удалось получить). Все это порождает довольно много разных сочетаний (вариантов) и требует большого числа оговорок. Этому планируется посвятить отдельную работу. Далее мы ограничимся одним специальном классом задач, играющих важную роль в моделировании компьютерных и транспортных сетей (см., например, [5], [35]).

**Пример 4.** Рассмотрим задачу ($Q$ – множество просто структуры, скажем, неотрицательный ортант)

$$\sum_{k=1}^{m} f_k\left(A_k^T x\right) + g(x) \to \min_{x \in Q},$$

где $g(x) = \sum_{i=1}^{n} g_i(x_i)$ (впрочем, часть изложенных далее конструкций не требует выполнения этого условия). Градиенты функции $f_k$ вычислимы за $\mathrm{O}(1)$ и имеют равномерно ограниченные константы (числом $L$) Липшица производной в 2-норме. Функция $g(x)$ предполагается сильно выпуклой в $p$-норме с константой $\mu_p$. Вводя матрицу $A = [A_1, \ldots, A_m]^T$ и вспомогательный вектор $z = Ax$ мы можем переписать эту задачу в "раздутом" пространстве $x := (x, z)$, как задачу типа проектирования на аффинное многообразие [35], рассмотренную ранее.[21] Однако для полноты картины[22] нам представляется полезнее провести для этой задачи рассуждения немного в другом ключе (следуя, например, [5], [35]). Прежде всего, заметим, что в эту схему погружаются следующие задачи [5]:

1) $\dfrac{L}{2}\|Ax - b\|_2^2 + \dfrac{\mu}{2}\|x - x_g\|_2^2 \to \min_{x \in \mathbb{R}^n}$,

2) $\dfrac{L}{2}\|Ax - b\|_2^2 + \mu \sum_{k=1}^{n} x_k \ln x_k \to \min_{x \in S_n(1)}$.

Константа Липшица производных одинаковы $L^1 = L^2 = L$, константы сильной выпуклости (считаются в разных нормах) также одинаковы $\mu_2^1 = \mu_1^2 = \mu$. Опишем далее довольно общий способ построения двойственной задачи:[23]

---

[21] В связи с этим, можно добавить к ограничениям, например, такого типа неравенства $A_k^T x \geq c_k$. Сложность задачи это не изменит. Этот факт можно использовать при численном поиске стохастических равновесий в модели стабильной динамики [35].

[22] Приводимая далее конструкция позволяет (с помощью перехода к двойственной задаче и ее последующего изучения покомпонентными методами) в некотором смысле перейти от игры на разной гладкости по разным направлениям для исходной задачи к и игре на разной сильной выпуклости функционала исходной задачи по разным направлениям (при переходе к двойственной задаче эта игра переходит в игру на гладкости двойственного функционала, которая уже неплохо проработана покомпонентными методами). К сожалению, все это возможно не в общем случае. Более того, в данной работе мы не будем подробно погружаться в детали. Этому планируется посвятить отдельную работу.

[23] Здесь важна оговорка о возможности "эффективно" решать задачу вида

$$\langle c, x \rangle + g(x) \to \max_{x \in Q}.$$

Вообще говоря, оговорка нетривиальная. В общем случае эта задача по сложности может соответствовать исходной. Стоит, однако, оговориться, что $g(x)$ в таких постановках в типичных приложениях является, как правило, "регуляризатором" исходной задачи (введенном нами с целью получения сильно выпуклой постановки или, например, возникшем при байесовском оценивании в качестве прайера или просто как пе-

$$\min_{x \in Q} \left\{ \sum_{k=1}^{m} f_k \left( A_k^T x \right) + g(x) \right\} = \min_{\substack{x \in Q \\ z = Ax}} \left\{ \sum_{k=1}^{m} f_k \left( z_k \right) + g(x) \right\} =$$

$$= \min_{\substack{x \in Q \\ z = Ax, z'}} \max_{y} \left\{ \langle z - z', y \rangle + \sum_{k=1}^{m} f_k \left( z_k' \right) + g(x) \right\} =$$

$$= \max_{y \in \mathbb{R}^m} \left\{ -\max_{\substack{x \in Q \\ z = Ax}} \left\{ \langle -z, y \rangle - g(x) \right\} - \max_{z'} \left\{ \langle z', y \rangle - \sum_{k=1}^{m} f_k \left( z_k' \right) \right\} \right\} =$$

$$= \max_{y \in \mathbb{R}^m} \left\{ -\max_{x \in Q} \left( \langle -A^T y, x \rangle - g(x) \right) - \sum_{k=1}^{m} \max_{z_k'} \left( z_k' y_k - f_k \left( z_k' \right) \right) \right\} =$$

$$= \max_{y \in \mathbb{R}^m} \left\{ -g^* \left( -A^T y \right) - \sum_{k=1}^{m} f_k^* \left( y_k \right) \right\} = -\min_{y \in \mathbb{R}^m} \left\{ g^* \left( -A^T y \right) + \sum_{k=1}^{m} f_k^* \left( y_k \right) \right\}.$$

Для упомянутых задач, получим:

1) $\dfrac{1}{2\mu} \left( \left\| x_g - A^T y \right\|_2^2 - \left\| x_g \right\|_2^2 \right) + \dfrac{1}{2L} \left( \left\| y + b \right\|_2^2 - \left\| b \right\|_2^2 \right) \to \min_{y \in \mathbb{R}^m}$,

2) $\dfrac{1}{\mu} \ln \left( \sum_{i=1}^{n} \exp \left( \dfrac{\left[ -A^T y \right]_i}{\mu} \right) \right) + \dfrac{1}{2L} \left( \left\| y + b \right\|_2^2 - \left\| b \right\|_2^2 \right) \to \min_{y \in \mathbb{R}^m}$.

В общем случае можно утверждать, что $\sum_{k=1}^{m} f_k^* \left( y_k \right)$ (композитный член в двойственной задаче) является сильно выпуклым в стандартной евклидовой норме (2-норме) с константой сильной выпуклости равной $L^{-1}$. Легко понять, что изучение свойств гладкости $g^* \left( -A^T y \right)$ (с точностью до множителя $\mu^{-1}$) совершенно аналогично тому, что мы уже делали в примере 3. То есть можно утверждать, что для двойственной задачи в стандартной евклидовой норме (2-норме)

$$\breve{L}_{\text{БГМ}} = \frac{1}{\mu} \max_{\|y\|_2 \le 1, \|x\|_p \le 1} \langle A^T y, x \rangle^2 = \frac{1}{\mu} \max_{\|x\|_p \le 1} \|Ax\|_2^2 = \frac{1}{\mu} \begin{cases} 1) \; \lambda_{\max} \left( A^T A \right) \\ 2) \; \max_{k=1,\ldots,n} \left\| A^{\langle k \rangle} \right\|_2^2 \end{cases},$$

---

нализация за сложность модели и т.д.). В любом случае, мы, как правило, имеем достаточно степеней свободы, чтобы добиться нужной простоты этой вспомогательной задачи. Сепарабельность $g(x)$ здесь довольно часто является "ключом к успеху".

$$\breve{L}_{\textit{ПБГМ}} \le \frac{1}{\mu} \max_{\|y\|_1 \le 1, \|x\|_p \le 1} \left\langle A^T y, x \right\rangle^2 = \frac{1}{\mu} \max_{\|x\|_p \le 1} \|Ax\|_\infty^2 = \frac{1}{\mu} \begin{cases} 1) \ \max_{k=1,\ldots,m} \|A_k\|_2^2 \\ 2) \ \max_{\substack{i=1,\ldots,m \\ j=1,\ldots,n}} |A_{ij}|^2 \end{cases}.$$

Для метода ACRCD* из замечания 2 (если смотреть через замечание 6, то следует полагать $\beta = 0$), примененного к двойственной задаче в не разреженном случае имеем следующие оценки времени работы:

1) $T_1 = \tilde{O}\left( nm \sqrt{\dfrac{L \max_{k=1,\ldots,m} \|A_k\|_2^2}{\mu}} \right)$,

2) $T_2 = \tilde{O}\left( nm \sqrt{\dfrac{L \max_{i,j} |A_{ij}|^2}{\mu}} \right)$.

Если теперь посмотреть на исходную прямую задачу (с $L := L/m$)

$$\frac{1}{m} \sum_{k=1}^{m} f_k\left(A_k^T x\right) + g(x) \to \min_{x \in Q},$$

и оценить время работы ускоренного метода рандомизации суммы из п. 3 данной статьи, то получим анонсированное в п. 3 соответствие с приведенными только что оценками (с $L := L/m$). Действительно, с учетом того, что константа Липшица градиентов $f_k\left(A_k^T x\right)$, посчитанные в соответствующих нормах (соответствующей норме, в которой сильно выпукл композит прямой задачи), равномерно оцениваются следующим образом:

1) $L \max_{k=1,\ldots,m} \|A_k\|_2^2$,

2) $L \max_{i,j} |A_{ij}|^2$,

а сложность вычисления $\nabla f_k\left(A_k^T x\right)$ равна $O(n)$, то согласно п. 3 имеем следующие оценки времени работы (далее считаем, что $m$ меньше $\min\{\ \}$):[24]

---

[24] Можно показать, развивая конструкцию п. 6.3 [1] с помощью идей, изложенных в замечании 8, что неускоренная составляющая этой оценки допускает для определенного класса задач в разреженном случае стоимости итерации $\tilde{O}(s)$.

$$1)\quad \tilde{T}_1 = \tilde{O}\left( n \cdot \left( m + \min\left\{ \sqrt{m \frac{L \max_{k=1,..,m} \|A_k\|_2^2}{\mu}}, \frac{L \max_{k=1,..,m} \|A_k\|_2^2}{\mu} \right\} \right) \right) = \tilde{O}\left( nm\sqrt{\frac{(L/m) \max_{k=1,..,m} \|A_k\|_2^2}{\mu}} \right),$$

$$2)\quad \tilde{T}_2 = \tilde{O}\left( n \cdot \left( m + \min\left\{ \sqrt{m \frac{L \max_{i,j} |A_{ij}|^2}{\mu}}, \frac{L \max_{i,j} |A_{ij}|^2}{\mu} \right\} \right) \right) = \tilde{O}\left( nm\sqrt{\frac{(L/m) \max_{i,j} |A_{ij}|^2}{\mu}} \right).$$

Таким образом, имеет место полное соответствие (с точностью до опущенных в рассуждениях логарифмических множителей).[25] Интересно заметить, что для первой задачи здесь также как и в примере 3 можно сполна использовать разреженность матрицы $A$. Более того, эту задачу (также с полным учетом разреженности) можно решать и прямым ПБГМ (см. замечание 8). Соответствующие оценки имеют вид [5] (мы снова возвращаемся к исходному пониманию параметра $L$):

$$T_1^{прям} = \tilde{O}\left( sn\sqrt{\frac{L \max_{k=1,...,n} \|A^{\langle k \rangle}\|_2^2}{\mu}} \right),$$

$$T_1^{двойств} = \tilde{O}\left( \tilde{s}m\sqrt{\frac{L \max_{k=1,..,m} \|A_k\|_2^2}{\mu}} \right) = \tilde{O}\left( sn\sqrt{\frac{L \max_{k=1,..,m} \|A_k\|_2^2}{\mu}} \right).$$

В действительности, обе эти оценки оказываются завышенными.[26] Более аккуратные рассуждения, позволяют обобщить результаты примеров 3, 4 на случай неравноправия слагаемых. Все это приведет к замене максимума на некоторое (в зависимости от выбора $\beta$) среднее. Скажем, в упомянутых уже ранее транспортных приложениях [5], [32], [37]–[35] матрица $A$ не просто разреженная, но еще и битовая (состоит из нулей и единиц). В таком случае приведенные оценки переписываются следующим образом:

---

[25] Соответствие имеет место и для неускоренной составляющей выписанных оценок (собственно, именно этот случай рассматривался в п. 3 данной статьи). Чтобы это понять, нужно в оценки для неускоренных покомпонентных методов $T_{ПМ}$ (см. текст сразу после замечания 2 и пример 3 с $s = m$, $\tilde{s} = n$)

$$T_{ПМ}^1 = \tilde{O}\left( nm \frac{L \max_{k=1,..,m} \|A_k\|_2^2}{\mu} \right),\quad T_{ПМ}^2 = \tilde{O}\left( nm \frac{L \max_{i,j} |A_{ij}|^2}{\mu} \right)$$

подставить $L := L/m$.

[26] Это легко усмотреть из способа рассуждений, в котором мы заменяем константы Липшица по разным направлениям, на худшую из них, это также позволяет эффективно использовать оценку ACRCD* с $\beta = 0$, см. замечания 2, 6.

$$T_1^{прям} = \tilde{O}\left(sn\sqrt{\frac{Ls}{\mu}}\right), \quad T_1^{двойств} = \tilde{O}\left(sn\sqrt{\frac{L\tilde{s}}{\mu}}\right).$$

Отсюда можно сделать довольно неожиданный вывод [5]: при $m \ll n$ стоит использовать прямой ПБГМ, а в случае $m \gg n$ – двойственный. Первый случай соответствует приложениям к изучению больших сетей (компьютерных, транспортных). Второй случай соответствует задачам, приходящим из анализа данных. □

### 6. Заключительные замечания

Если сравнить общие трудозатраты быстрого градиентного метода и его покомпонентного варианта, то довольно часто удается ускориться в $\sim \sqrt{n}$ раз (где $n$ – размерность пространства, в котором происходит оптимизация). Собственно, значительная часть данной работы (пп. 4, 5) была посвящена обсуждению того, в каких ситуациях можно рассчитывать на такое ускорение.

Отмеченное ускорение происходит за счет "обмана" потенциального сопротивляющегося оракула, корректирующего исходя из оставшихся у него свобод по ходу итерационного процесса оптимизируемую функцию таким образом, чтобы нам доставались наиболее плохие направления градиента (с большой константой Липшица – это соответствует пилообразному движению по дну растянутого оврага, с медленным приближением к середине оврага, в которой расположился минимум).

Введение рандомизации в метод – это универсальный рецепт гарантированно обезопасить себя от худшего случая. Причем важно отметить, что это не теоретический трюк, который позволяет просто гарантировать лучшую теоретическую оценку. К сожалению, овражность функций – это довольно типичное свойство задач больших размеров, поэтому даже если мы возьмем сложность в среднем (по множеству типичных входов) для быстрого градиентного метода, то оценка получится все хуже (поскольку типичные входы не столь хороши), чем для рандомизированного метода. В обоих случаях мы считаем средние затраты (математическое ожидание), только в разных пространствах и по разным вероятностным мерам. В случае рандомизированного метода мы частично диверсифицируем себя от всего того, что может быть на входе, но это происходит не бесплатно, а с помощью препроцессинга, требующего расчет констант Липшица градиента по всем направлениям. К счастью, такой препроцессинг можно делать адаптивно и эффективно.[27]

---

[27] Это можно делать за небольшую дополнительную плату – мультипликативный фактор порядка 4 в числе обращений за $\partial f(x)/\partial x_i$, если рассматривается композитная постановка или $Q$ не параллелепипедного типа, то требуется еще уметь рассчитывать (пересчитывать) вместе с $\partial f(x)/\partial x_i$ и значение функции.

Авторы выражают благодарность Ю.Е. Нестерову, Александру Рахлину и Питеру Рихтарику за ряд ценных ссылок, а также Александру Тюрину и Дмитрию Камзолову за помощь в работе.

Настоящая статья представляет собой запись совместного доклада А.В. Гасникова и И.Н. Усмановой на "International Conference on Operation Research 2015. Optimal Decision and Big Data." Vienna, September 1 – 4, 2015.

В конце марта 2016 (когда данная статья уже ожидала своей очереди публикации в журнале Труды МФТИ) авторы узнали, что параллельно с настоящей статьей появился электронный препринт [41], в котором получены довольно близкие результаты к части результатов, приведенных в данной статье (замечание 6 при $\beta = 1/2$ и пример 2). Однако для авторов этой работы первоисточником идеи (для замечания 6 при $\beta = 1/2$ и примера 2) и отправной точкой в развитие этого направления послужило выступление Ю.Е. Нестерова [27]. Заметим также, что в марте 2016 появился препринт [42], в основу которого положено выступление [27].

В феврале 2016 года появился электронный препринт [43], в котором обсуждаются общие прямо-двойственные подходы для постановки задачи из примера 4.